\documentclass[12pt,a4paper]{article}

\usepackage{amsmath}
\usepackage{amsthm}
\usepackage{amssymb}
\usepackage{amscd}
\usepackage{tikz}

\usepackage{a4wide}

\theoremstyle{definition}
\newtheorem{theorem}{Theorem}[section]
\newtheorem{prop}[theorem]{Proposition}
\newtheorem{lemma}[theorem]{Lemma}
\newtheorem{corollary}[theorem]{Corollary}
\newtheorem{definition}[theorem]{Definition}
\newtheorem{example}[theorem]{Example}
\newtheorem{remark}[theorem]{Remark}

\newenvironment{demo}[1]{%
  \trivlist
  \item[\hskip\labelsep
        {\bf #1.}]
}{%
  \endtrivlist
}

\renewcommand\tilde{\widetilde}

\newcommand\ep{\varepsilon}

\newcommand\AC{\mathcal{A}}
\newcommand\CC{\mathcal{C}}
\newcommand\EE{\mathcal{E}}
\newcommand\FF{\mathcal{F}}
\newcommand\OO{\mathcal{O}}
\newcommand\PP{\mathcal{P}}
\newcommand\ST{\mathcal{S}}
\newcommand\TT{\mathcal{T}}
\newcommand\eAC{\mathcal{A}^{(e)}} 
\newcommand\eCC{\mathcal{C}^{(e)}}
\newcommand\eFF{\mathcal{F}^{(e)}}

\newcommand\eOO{\mathcal{O}^{(e)}}
\newcommand\eST{\mathcal{S}^{(e)}}
\newcommand\eTT{\mathcal{T}^{(e)}}
\newcommand\ePhi{\Phi^{(e)}}
\newcommand\ePsi{\Psi^{(e)}}

\newcommand\K{\mathbb{K}}
\newcommand\vectx{\boldsymbol{x}}

\newcommand\Real{\mathbb{R}}
\newcommand\Int{\mathbb{Z}}
\newcommand\conv{\operatorname{conv}} 
\newcommand\supp{\operatorname{supp}} 
\newcommand\SC{\operatorname{SC}}
\newcommand\MC{\operatorname{MC}}
\newcommand\C{\operatorname{C}}
\newcommand\sgn{\operatorname{sgn}}
\newcommand\ctop{\operatorname{top}}

\allowdisplaybreaks

\title{
Two enriched poset polytopes
}

\author{
Soichi Okada\footnote{
Graduate School of Mathematics, Nagoya University, 
Furo-cho, Chikusa-ku, Nagoya 464-8602, Japan, 
{\tt okada@math.nagoya-u.ac.jp}
}
\ and
Akiyoshi Tsuchiya\footnote{
Department of Information Science,
	Faculty of Science,
	Toho University,
	Miyama, Funabashi-shi, Chiba 274-8510, Japan,
{\tt akiyoshi@is.sci.toho-u.ac.jp}
}
}

\date{}

\begin{document}

\maketitle

\abstract{%
Stanley introduced and studied two lattice polytopes, the order polytope and chain polytope, 
associated to a finite poset.
Recently Ohsugi and Tsuchiya introduce an enriched version of them, 
called the enriched order polytope and enriched chain polytope.
In this paper, we give a piecewise-linear bijection between these enriched poset polytopes, 
which is an enriched analogue of Stanley's transfer map 
and bijectively proves that they have the same Ehrhart polynomials.
Also we construct explicitly unimodular triangulations of two enriched poset polytopes, 
which are the order complexes of graded posets. 
\par
Mathematics Subject Classification (MSC2010): 52B12 (primary), 05E45, 06A07, 52B20 (secondary)
}

\section{%
Introduction
}

We assume that readers are familiar with the definition of a poset presented in \cite[Chapter 3]{EC1}.
Let $P$ be a finite poset with $d$ elements.
We denote by $\Real^P$ the vector space of all real-valued functions on $P$, 
and identify $\Real^P$ with the Euclidean space $\Real^d$.
The \emph{order polytope} $\OO(P)$ of $P$ is the subset of $\Real^P$ 
consisting of all functions $f : P \to \Real$ satisfying the following two conditions:
\begin{enumerate}
\item[(i)]
$0 \le f(v) \le 1$ for all $v \in P$;
\item[(ii)]
If $x < y$ in $P$, then we have $f(x) \le f(y)$.
\end{enumerate}
And the \emph{chain polytope} $\CC(P)$ of $P$ is the subset of $\Real^P$ consisting of 
all functions $g : P \to \Real$ satisfying the following two conditions:
\begin{enumerate}
\item[(i)]
$g(v) \ge 0$ for all $v \in P$;
\item[(ii)]
If $v_1 > \dots > v_r$ is a chain in $P$, then we have $g(v_1) + \dots + g(v_r) \le 1$.
\end{enumerate}
Then it is known (see \cite[Corollary~1.3 and Theorem~2.2]{Stanley1986}) that 
$\OO(P)$ and $\CC(P)$ are convex polytopes whose vertex sets are given by 
\begin{align*}
\FF(P)
 &=
\{ \chi_{F} : \text{$F$ is an order filter of $P$} \},
\\
\AC(P)
 &=
\{ \chi_A : \text{$A$ is an antichain of $P$} \}
\end{align*}
respectively, where $\chi_S$ is the characteristic function of a subset $S \subset P$ 
defined by $\chi_S(v) = 1$ if $v \in S$ and $0$ otherwise.
Here 
an \emph{order filter} of $P$ is a subset $F \subset P$ 
such that if $v \in F$ and $v < w$ then $w \in F$.
In particular, we have
\[
\OO(P) = \conv \FF(P),
\quad
\CC(P) = \conv \AC(P),
\]
where $\conv S$ denotes the convex hull of $S$.
These poset polytopes are related via the transfer map.

\begin{theorem}
\label{thm:transfer}
(Stanley \cite[Theorem~3.2]{Stanley1986})
We define a piecewise-linear map $\Phi : \Real^P \to \Real^P$, called the \emph{transfer map}, by
\begin{equation}
\label{eq:def_transfer}
\left( \Phi f \right)(v)
 =
\begin{cases}
f(v) &\text{if $v$ is minimal in $P$,} \\
f(v) - \max \{ f(w) : \text{$v$ covers $w$ in $P$} \}
 &\text{if $v$ is not minimal in $P$}
\end{cases}
\end{equation}
for $f \in \Real^P$ and $v \in P$.
Then $\Phi$ induces a continuous bijection from $\OO(P)$ to $\CC(P)$.
In particular, $\Phi$ provides a bijection between $m \OO(P) \cap \Int^P$ 
and $m \CC(P) \cap \Int^P$ for any nonnegative integer $m$, 
where $m \PP = \{ m f : f \in \PP \}$ is the $m$th dilation of a polytope $\PP$ 
and $\Int^P$ is the set of all integer-valued functions on $P$.
\end{theorem}

The transfer map enables us to compare certain properties of $\OO(P)$ and $\CC(P)$.
For example, the two polytopes $\OO(P)$ and $\CC(P)$ have the same Ehrhart polynomials, i.e.,
\begin{equation}
\label{eq:sameEhrhart}
\# \left( m \OO(P) \cap \Int^P \right)
 =
\# \left( m \CC(P) \cap \Int^P \right).
\end{equation}
We note that the polynomial $\# \left( m \OO(P) \cap \Int^P \right)$ in $m$ 
is the order polynomial (with a shifted argument) of the poset $P$, 
which counts the number of $P$-partitions. 
A map $h : P \to \Int_{\geq 0}$, where $\Int_{\geq 0}$ is the set of nonnegative integers, 
is called a \emph{$P$-partition} if $v \leq w$ implies $h(v) \leq h(w)$.
Then $m \OO(P) \cap \Int^P$ is the set of all $P$-partitions $h : P \to \Int_{\ge 0}$ 
such that $h(v) \le m$ for all $v \in P$.
Hence the transfer map $\Phi$ gives a bijection between such $P$-partitions and 
lattice points in the $m$th dilation of the chain polytope $\CC(P)$.
In a very recent work, Higashitani \cite{Higashitani2020} proves 
that $\OO(P)$ and $\CC(P)$ are combinatorially mutation-equivalent 
by using the transfer map $\Phi$. 
The notion of combinatorial mutation was introduced from viewpoints 
of mirror symmetry for Fano manifolds. 
Also we can transfer a canonical triangulation of $\OO(P)$, which is the order complex of a graded poset as simplicial complexes, to $\CC(P)$ via the transfer map $\Phi$.

\begin{theorem}
\label{thm:triangulation}
(Stanley \cite[Section~5]{Stanley1986})
For a chain $C = \{ F_1 \supsetneq F_2 \supsetneq \dots \supsetneq F_k \}$ 
of order filters of $P$, we put
\begin{equation}
\label{eq:simplex}
S_C = \conv \{ \chi_{F_1}, \dots, \chi_{F_k} \},
\quad
T_C = \conv \{ \Phi(\chi_{F_1}), \dots, \Phi(\chi_{F_k}) \}.
\end{equation}
Then we have
\begin{enumerate}
\item[(a)]
The collection $\ST_P = \{ S_C : \text{$C$ is a chain of order filters} \}$ is a unimodular triangulation 
of the order polytope $\OO(P)$.
\item[(b)]
The collection $\TT_P= \{ T_C : \text{$C$ is a chain of order filters} \}$ is a unimodular triangulation 
of the chain polytope $\CC(P)$.
\end{enumerate}
\end{theorem}

In the last decade, many authors have generalized order polytopes and chain polytopes 
together with their connecting transfer maps.
These generalizations include marked poset polytopes \cite{ABS2011} 
and double poset polytopes \cite{CFS2017, HT2017}.

Recently from a viewpoint of (left) enriched $P$-partitions, Ohsugi--Tsuchiya \cite{OT1, OT2} introduced an enriched version of 
order polytopes and chain polytopes as follows.
We define $\eFF(P)$ and $\eAC(P)$ by putting
\begin{align}
\label{eq:def_eF}
\eFF(P)
 &=
\left\{
 f \in \Int^P
 : 
\begin{array}{l}
 \text{(i) $f(v) \in \{ 1, 0, -1 \}$ for any $v \in P$, } \\
 \text{(ii) $\supp (f)$ is an order filter of $P$, and} \\
 \text{(iii) if $f(v) = -1$, then $v$ is minimal in $\supp(f)$}
\end{array}
\right\},
\\
\label{eq:def_eA}
\eAC(P)
 &=
\left\{
 f \in \Int^P
 :
\begin{array}{l}
 \text{(i) $f(v) \in \{ 1, 0, -1 \}$ for any $v \in P$, and} \\
 \text{(ii) $\supp(f)$ is an antichain of $P$}
\end{array}
\right\},
\end{align}
where $\supp (f) = \{ v \in P : f(v) \neq 0 \}$.
Then the \emph{enriched order polytope} $\eOO(P)$ 
and the \emph{enriched chain polytope} $\eCC(P)$ are defined as 
the convex hulls of $\eFF(P)$ and $\eAC(P)$ respectively:
\[
\eOO(P) = \conv \eFF(P),
\quad
\eCC(P) = \conv \eAC(P).
\]
(Note that enriched chain polytopes were independently introduced 
by Kohl, Olsen and Sanyal \cite[Section~7]{KOS} under the name 
of \emph{unconditional chain polytopes}.)
Since $\FF(P) = \eFF(P) \cap \{ 0, 1 \}^P$ and $\AC(P) = \eAC(P) \cap \{ 0, 1 \}^P$, 
we have $\OO(P) \subset \eOO(P)$ and $\CC(P) \subset \eCC(P)$.
Then, 
as an enriched version of (\ref{eq:sameEhrhart}), 
Ohsugi--Tsuchiya \cite{OT2} used a commutative algebra technique to prove
\begin{equation}
\label{eq:e-sameEhrhart}
\# \left( m \eOO(P) \cap \Int^P \right)
 =
\# \left( m \eCC(P) \cap \Int^P \right)
\end{equation}
for any nonnegative integer $m$. 
It is a natural problem to find a bijective proof of this equality (\ref{eq:e-sameEhrhart}).

For a nonnegative integer $m$, we denote by $\EE_m(P)$ the set of all left enriched $P$-partitions 
$h : P \to \Int$ such that $|h(v)| \leq m$ for any $v \in P$
(see Section~\ref{subsec:enrchiedbijection} for the definition of left enriched $P$-partitions).
Then we have $\eFF(P) = \eOO(P) \cap \Int^P = \EE_1(P)$, 
and if $m \geq 2$, then a map $h \in m \eOO(P) \cap \Int^P$ is not always 
a left enriched $P$-partition (see \cite[Example 4.2]{OT2}). 
However it is known (\cite[Theorem~0.2 and its proof]{OT1}) that 
there is an explicit bijection between $m \eCC(P) \cap \Int^P$ and $\EE_m(P)$. This and (\ref{eq:e-sameEhrhart}) are the reasons why we call $\eOO(P)$ and $\eCC(P)$ ``enriched" poset polytopes.

One of the main results of this paper is the following theorem, 
which gives a bijective proof of (\ref{eq:e-sameEhrhart}).

\begin{theorem}
\label{thm:e-transfer}
We define a piecewise-linear map $\ePhi : \Real^P \to \Real^P$, which we call 
the \emph{enriched transfer map}, 
inductively on the ordering of $P$ such that the value of $\ePhi(f)$ at $v$ is equal to 
\begin{multline}
\label{eq:def_e-transfer}
\left\{
 \begin{array}{ll}
 f(v) &\text{if $v$ is minimal in $P$,}
 \\
 \multicolumn{2}{l}{
 f(v)
  - 
 \max \left\{ 
  \displaystyle\sum_{i=1}^r \left| \left( \ePhi f \right)(v_i) \right|
  : \text{$v > v_1 > \dots > v_r$ is a chain in $P$} \right\}
 }
 \\
 &\text{if $v$ is not minimal in $P$.}
 \end{array}
\right.
\end{multline}
Then $\ePhi$ induces a continuous bijection from $\eOO(P)$ to $\eCC(P)$.
In particular, $\ePhi$ provides a bijection between $m \eOO(P) \cap \Int^P$ and $m \eCC(P) \cap \Int^P$ 
for any nonnegative integer $m$.
\end{theorem}

Moreover, by composing with $\ePhi$, we also obtain an explicit bijection 
between $m \eOO(P) \cap \Int^P$ and $\EE_m(P)$ for any nonnegative integer $m$ 
(Proposition~\ref{prop:bij_P2O}). 

It can be shown (see Proposition~\ref{prop:e-transfer-O2C}) 
that the restriction of $\ePhi$ to $\OO(P)$ 
gives a continuous piecewise-linear bijection between $\OO(P)$ and $\CC(P)$, 
which coincides with the restriction of Stanley's transfer map $\Phi$ in Theorem~\ref{thm:transfer}.
Also, by the same technique of \cite{Higashitani2020}, 
we can show that $\eOO(P)$ and $\eCC(P)$ are combinatorially mutation-equivalent 
by using the enriched transfer map $\ePhi$ (see \cite[Section~5]{Higashitani2020}).

Ohsugi--Tsuchiya \cite{OT1,OT2} constructed triangulations of enriched order and chain polytopes 
by using the algebraic technique of Gr\"{o}bner bases.
Also, Kohl--Olsen--Sanyal \cite{KOS} constructed triangulations of enriched chain polytopes 
from a viewpoint of convex geometry.
Another main result of this paper is an explicit combinatorial description 
of triangulations of two enriched poset polytopes, 
which are the order complexes of graded posets as simplicial complexes and 
are transferred by the enriched transfer map $\ePhi$.
Our result is analogous to Stanley's canonical triangulations of two poset polytopes 
(see Theorem~\ref{thm:triangulation}).

\begin{theorem}
\label{thm:e-triangulation}
We equip $\eFF(P)$ with a poset structure by the partial ordering given in Definition~\ref{def:ordering}.
For a chain $K$ in $\eFF(P)$, we define
\begin{equation}
\label{eq:e-simplex}
S^{(e)}_K = \conv K, 
\quad
T^{(e)}_K = \conv \ePhi(K).
\end{equation}
Then we have
\begin{enumerate}
\item[(a)]
The set $\eST_P = \{ S^{(e)}_K : \text{$K$ is a chain in $\eFF(P)$} \}$ is a unimodular triangulation of $\eOO(P)$.
\item[(b)]
The set $\eTT_P = \{ T^{(e)}_K : \text{$K$ is a chain in $\eFF(P)$} \}$ is a unimodular triangulation of $\eCC(P)$.
\end{enumerate}
\end{theorem}

Remark that the partial ordering on $\eFF(P)$ given in Definition~\ref{def:ordering} 
is an extension of the inclusion ordering on the set of order filters of $P$, 
so the poset $\FF(P)$ is the induced subposet of $\eFF(P)$.
Stanley gave the defining inequalities of facets of the canonical triangulations $\ST_P$ and $\TT_P$ of $\OO(P)$ and $\CC(P)$ (\cite[Section~5]{Stanley1986}). We also give sets of defining inequalities of facets of the triangulation $\eST_P$ and $\eTT_P$ of $\eOO(P)$ and $\eCC(P)$ (Corollary~\ref{cor:e-facet-C} and Proposition~\ref{prop:e-facet-O}).
On the other hand, we identify these triangulations with Ohsugi--Tsuchiya's triangulations algebraically obtained in \cite{OT1,OT2} (Propositions \ref{prop:sametriangulation_eO} and \ref{prop:sametriangulation_eC}).

The rest of this paper is organized as follows.
In section~\ref{sec:transfer}, we prove Theorem~\ref{thm:e-transfer}, 
and give an explicit bijection between left enriched $P$-partitions and lattice points 
of the dilated enriched order polytope.
Section~\ref{sec:triangulation} is devoted to the proof of Theorem~\ref{thm:e-triangulation}.
We also give sets of defining inequalities for the maximal faces.
In Section~\ref{sec:algebraic}, we prove that the triangulations described in Theorem~\ref{thm:e-triangulation} 
coincide with the Ohsugi--Tsuchiya's triangulations.

\subsection*{Acknowledgements}
The first author was partially supported by JSPS KAKENHI 18K03208 
and the second author was partially supported by JSPS KAKENHI 19J00312, 19K14505 and 22K13890.

\section{%
Enriched transfer map
}
\label{sec:transfer}

In this section, we give a proof of Theorem~\ref{thm:e-transfer}, 
and we use the enriched transfer map to describe a bijection between 
left enriched $P$-partitions and lattice points of the dilated enriched order polytope.
 
\subsection{%
Notations
}

In what follows, we use the following notations and terminologies.
Let $P$ be a finite poset.
For $v$, $w \in P$, we say that $v$ covers $w$, written $v \gtrdot w$, 
if $v > w$ and there is no element $u$ such that $v > u > w$.
Given an antichain $A$, we denote by $\langle A \rangle$ 
the smallest order filter containing $A$.
Given an element $v \in P$, we put
\[
P_{\le v} = \{ w \in P : w \le v \},
\quad
P_{<v} = \{ w \in P : w < v \}.
\]
For a subposet $Q$ of $P$, 
we denote by $\max Q$ and $\min Q$ the set of maximal and minimal elements of $Q$ respectively.
For a chain $C = \{ v_1 > v_2 > \dots > v_r \}$ of $Q$, we say that
\begin{itemize}
\item
$C$ is \emph{saturated} if $v_i \gtrdot v_{i+1}$ for $i=1, \dots, r-1$;
\item
$C$ is \emph{maximal} if it is saturated and $v_1 \in \max Q$ and $v_r \in \min Q$.
\end{itemize}
Let $\C(Q)$, $\SC(Q)$ and $\MC(Q)$ be 
the sets of all chains, all saturated chains and all maximal chains 
respectively.
We denote by $\ctop C$ the maximum element of a chain $C$.
For $f \in \Real^P$ and a chain $C = \{ v_1 > \dots > v_r \}$, we define
\begin{align*}
S(f;C)
 &=
|f(v_1)| + \dots + |f(v_r)|,
\\
T^+(f;C)
 &=
- f(v_1) - 2 f(v_2) - \cdots - 2^{r-2} f(v_{r-1}) + 2^{r-1} f(v_r),
\\
T^-(f;C)
 &=
- f(v_1) - 2 f(v_2) - \cdots - 2^{r-2} f(v_{r-1}) - 2^{r-1} f(v_r).
\end{align*}
Note that, if $C$ is a one-element chain $\{ v \}$, 
then $T^+(f ; \{ v \}) = f(v)$ and $T^-(f ; \{ v \}) = - f(v)$.

\subsection{%
Defining inequalities for enriched poset polytopes
}

Our proof of Theorem~\ref{thm:e-transfer} is based on the defining inequalities of $\eOO(P)$ and $\eCC(P)$ 
given by \cite{OT2}.

\begin{prop}
\label{prop:ineq}
(\cite[Lemma~1.1]{OT1}, \cite[Proposition~6.1 and Theorem~6.2]{OT2})
We have
\begin{multline}
\label{eq:ineq_eO}
\eOO(P)
\\
 =
\left\{
 f \in \Real^P :
 \begin{array}{l}
  \text{$T^+(f;C) \le 1$ for all $C \in \SC(P)$ with $\ctop C \in \max (P)$} \\
  \text{$T^-(f;C) \le 1$ for all $C \in \MC(P)$}
 \end{array}
\right\},
\end{multline}
and
\begin{equation}
\label{eq:ineq_eC}
\eCC(P)
 =
\left\{
 g \in \Real^P :
 \text{$S(g;C) \le 1$ for all $C \in \MC(P)$}
\right\}.
\end{equation}
\end{prop}

\begin{example}
Let $\Lambda$ be the three-element poset on $\{ u, v, w \}$ 
with covering relations $u \lessdot w$ and $v \lessdot w$.
If we identify $\Real^\Lambda$ with $\Real^3$ by the correspondence 
$f \leftrightarrow (f(u), f(v), f(w))$, we have
\begin{align*}
\eFF(\Lambda)
 &=
\left\{
\begin{array}{l}
 (0,0,0), (0,0,1), (0,0,-1), (1,0,1), (-1,0,1), (0,1,1), (0,-1,1) \\
 (1,1,1), (1,-1,1), (-1,1,1), (-1,-1,1)
\end{array}
\right\},
\\
\eAC(\Lambda)
 &=
\left\{
\begin{array}{l}
 (0,0,0), (1,0,0), (-1,0,0), (0,1,0), (0,-1,0), (0,0,1), (0,0,-1) \\
 (1,1,0), (1,-1,0), (-1,1,0), (-1,-1,0)
\end{array}
\right\},
\end{align*}
and
\begin{align*}
\eOO(\Lambda)
 &=
\left\{
 f \in \Real^\Lambda :
\begin{array}{l}
 f(w) \le 1 \\
 -f(u) + 2 f(w) \le 1, \ -f(v) + 2 f(w) \le 1 \\
 -f(u) - 2 f(w) \le 1, \ -f(u) - 2 f(w) \le 1
\end{array}
\right\},
\\
\eCC(\Lambda)
 &=
\left\{
 g \in \Real^\Lambda :
 |g(u)| + |g(w)| \le 1, \ |g(v)| + |g(w)| \le 1
\right\}.
\end{align*}
\end{example}

\subsection{%
Proof of Theorem~\ref{thm:e-transfer}
}

In this subsection we prove Theorem~\ref{thm:e-transfer}.
The inductive definition (\ref{eq:def_e-transfer}) of $\ePhi$ can be written as
\begin{multline}
\label{eq:e-transfer}
\left( \ePhi (f) \right)(v)
\\
 =
\begin{cases}
 f(v) &\text{if $v$ is minimal in $P$,} \\
 f(v) - \max\{ S(\ePhi(f);C) : C \in \C(P_{<v}) \} &\text{if $v$ is not minimal in $P$.}
\end{cases}
\end{multline}
It is easy to see that the map $\ePhi : \Real^P \to \Real^P$ is a bijection.

\begin{lemma}
\label{lem:inv_transfer}
The map $\ePhi : \Real^P \to \Real^P$ is a bijection with inverse map $\ePsi$ given by
\begin{multline}
\label{eq:e-invtransfer}
\left( \ePsi (g) \right)(v)
\\
 =
\begin{cases}
g(v) &\text{if $v$ is minimal in $P$,} \\
g(v)
 + 
\max \{ S(g ; C) :  C \in \C(P_{<v}) \}
 &\text{if $v$ is not minimal in $P$.}
\end{cases}
\end{multline}
\end{lemma}

Here we note that
\[
\max \{ S(g ; C) : C \in \C(P_{<v}) \}
 =
\max \{ S(g ; C) : C \in \MC(P_{<v}) \},
\]
hence we may replace $\C(P_{<v})$ with $\MC(P_{<v})$ in (\ref{eq:e-transfer}) 
and (\ref{eq:e-invtransfer}).
The following proposition follows from the definitions of $\ePhi$ and $\ePsi$.

\begin{prop}
\label{prop:e-transfer-F2A}
\begin{enumerate}
\item[(a)]
For $f \in \eFF(P)$, we have
\[
(\ePhi (f))(v)
 =
\begin{cases}
 f(v) &\text{if $v$ is minimal in $\supp(f)$,} \\
 0 &\text{otherwise.}
\end{cases}
\]
In particular, $\ePhi(f) \in \eAC(P)$ and $\supp \ePhi(f) = \min (\supp(f))$.
\item[(b)]
For $g \in \eAC(P)$, we have
\[
(\ePsi(g))(v)
 =
\begin{cases}
 1 &\text{if $v \in \langle \supp(g) \rangle \setminus \min \langle \supp(g) \rangle$,} \\
 g(v) &\text{if $v \in \min \langle \supp(g) \rangle$,} \\
 0 &\text{otherwise.}
\end{cases}
\]
In particular, $\ePsi(g) \in \eFF(P)$ and $\supp \ePsi(g) = \langle \supp(g) \rangle$.
\item[(c)]
The map $\ePhi$ induces a bijection between $\eFF(P)$ and $\eAC(P)$.
\end{enumerate}
\end{prop}

In order to prove Theorem~\ref{thm:e-transfer}, we need to prepare two lemmas.
We put
\begin{align*}
M(g;P_{\le v})
 &=
\max \{ S(g ; C) : C \in \MC(P_{\le v}) \},
\\
M(g;P_{< v})
 &=
\max \{ S(g ; C) : C \in \MC(P_{< v}) \}.
\end{align*}

\begin{lemma}
\label{lem:2-1}
Let $f \in \Real^P$ and $v \in P$.
We put
\[
\TT(f ; v)
 =
\{ T^+(f;C) : C \in \SC(P_{\le v}) \text{ with } \ctop C = v \}
 \cup
\{ T^-(f;C) : C \in \MC(P_{\le v}) \}.
\]
Then, for any $C \in \MC(P_{\le v})$,
there exists an element $T \in \TT(f ; v)$ such that $S(\ePhi(f); C) \le T$.
\end{lemma}

\begin{demo}{Proof}
We write $g = \ePhi(f)$.
We proceed by induction on the ordering of $P$.
If $v$ is a minimal element, then $C$ is a one-element chain $\{ v \}$ and
\[
S(g;C) = |g(v)| = |f(v)|
 =
\begin{cases}
 f(v) = T^+(f;C) &\text{if $f(v) \ge 0$,} \\
 -f(v) = T^-(f;C) &\text{if $f(v) \le 0$.}
\end{cases}
\]
If $v$ is not a minimal element, then by definition
\[
g(v) = f(v) - M(g;P_{<v}).
\]
Let $C = \{ v = v_1 \gtrdot v_2 \gtrdot \cdots \gtrdot v_r \}$.
Since $C \setminus \{ v \}  = \{ v_2 \gtrdot \cdots \gtrdot v_r \} \in \MC(P_{< v})$, 
we have
\[
S(g ; C \setminus \{ v \})
 \le 
M(g ; P_{<v}).
\]
If $g(v) = f(v) - M(v;P_{<v}) \ge 0$, then we have
\begin{align*}
S(g ; C)
 &=
f(v) - M(g; P_{<v})
 +
S(g ; C \setminus \{ v \})
 \\
&
 \le
f(v) = T^+(f ; \{ v \}).
\end{align*}
If $g(v) \le 0$, then we have
\begin{align*}
S(g ; C)
 &=
- f(v) + M(g; P_{<v})
 +
S(g ; C \setminus \{ v \})
\\
&
 \le
- f(v) 
+
2 M(g; P_{<v}).
\end{align*}
Let $C' \in \MC(P_{<v})$ be a chain which attains the maximum $M(g;P_{<v})$. 
Then, by applying the induction hypothesis to $C'$ and $w = \ctop C'$, 
there exists a chain $C''$ satisfying one of the following conditions:
\begin{enumerate}
\item[(i)]
$C'' \in \SC(P_{\le w})$ with $\ctop C'' = w$ and $S(g ; C') \le T^+(g ; C'')$;
\item[(ii)]
$C'' \in \MC(P_{\le w})$ and $S(g ; C') \le T^-(g ; C'')$.
\end{enumerate}
In the case (i), we have
\[
S(g ; C) \le - f(v) + 2 S(g ; C') \le - f(v) + 2 T^+(g ; C'') = T^+(g ; \{ v \} \cup C''),
\]
and in the case (ii), we have
\[
S(g ; C) \le - f(v) + 2 S(g ; C') \le - f(v) + 2 T^-(g ; C'') = T^-(g ; \{ v \} \cup C'').
\]
Since $v \gtrdot w$, we can complete the proof.
\qed
\end{demo}

\begin{lemma}
\label{lem:2-2}
Let $g \in \Real^P$ and $v \in P$.
For a chain $C = \{ v_1 \gtrdot v_2 \gtrdot \cdots \gtrdot v_r \} \in \SC(P_{\le v_1})$, 
we have
\[
2^{r-1} \left( |g(v_r)| + M(g; P_{< v_r}) \right)
+
\sum_{i=1}^{r-1}
 2^{r-i-1} \left( |g(v_{r-i})| - M(g; P_{< v_{r-i}}) \right)
\le
M(g; P_{\le v_1}).
\]
\end{lemma}

\begin{demo}{Proof}
We proceed by induction on $r$.
If $r = 1$, then
\begin{align*}
|g(v_1)| + M(g ; P_{< v_1})
 &=
|g(v_1)| + \max \{ S(g;C') : C' \in \MC(P_{< v_1}) \}
\\
 &=
\max \{ |g(v_1)| + S(g;C') : C' \in \MC(P_{< v_1}) \}
\\
 &=
\max \{ S(g;C) : C \in \MC(P_{\le v_1}) \}
\\
 &=
M(g;P_{\le v_1}).
\end{align*}

Let $r \ge 2$.
Since $\{ v_r \} \cup C' \in \MC(P_{< v_{r-1}})$ for any $C' \in \MC(P_{< v_r})$, 
we have
\begin{multline*}
|g(v_r)| + M(g; P_{< v_r})
 =
|g(v_r)| + \max \{ S(g;C') : C' \in MC(P_{< v_r}) \}
\\
 \le
\max \{ S(g;C'') : C'' \in \MC(P_{< v_{r-1}}) \}
 =
M(g;P_{< v_{r-1}}).
\end{multline*}
Hence we have
\begin{align*}
&
2^{r-1} \left( |g(v_r)| + M(g;P_{< v_r}) \right)
+
2^{r-2} \left( |g(v_{r-1})| - M(g; P_{< v_{r-1}}) \right)
\\
&
 \le
2^{r-1} M(g; P_{< v_{r-1}})
+
2^{r-2} \left( |g(v_{r-1})| - M(g; P_{< v_{r-1}}) \right)
\\
 &=
2^{r-2} \left( |g(v_{r-1})| + M(g; P_{< v_{r-1}}) \right).
\end{align*}
Therefore, by using the induction hypothesis, we see that
\begin{align*}
&
2^{r-1} \left( |g(v_r)| + M(g; P_{< v_r}) \right)
+
\sum_{i=1}^{r-1}
 2^{r-i-1} \left( |g(v_{r-i})| - M(g; P_{< v_{r-i}}) \right)
\\
&
 \le
2^{r-2} \left( |g(v_{r-1})| + M(g; P_{< v_{r-1}}) \right)
 +
\sum_{i=2}^{r-1}
 2^{r-i-1} \left( |g(v_{r-i})| - M(g; P_{< v_{r-i}}) \right)
\\
&
 \le
M(g; P_{\le v_1}).
\end{align*}
This completes the proof.
\qed
\end{demo}

Now we are ready to prove Theorem~\ref{thm:e-transfer}.

\begin{demo}{Proof of Theorem~\ref{thm:e-transfer}}
First we shall prove that $f \in \eOO(P)$ implies $\ePhi(f) \in \eCC(P)$.
Let $f \in \eOO(P)$ and put $g = \ePhi(f)$.
We show that $S(g;C) \le 1$ for all maximal chains $C = \{ v_1 \gtrdot v_2 \gtrdot \dots \gtrdot v_r \} \in \MC(P)$.
By Lemma~\ref{lem:2-1}, there exists a chain $C'$ satisfying one of the following conditions:
\begin{enumerate}
\item[(i)]
$C' \in \SC(P_{\le v_1})$ with $\ctop C' = v_1$ and $S(g;C) \le T^+(f;C')$;
\item[(ii)]
$C' \in \MC(P_{\le v_1})$ and $S(g;C) \le T^-(f;C')$.
\end{enumerate}
Then it follows from (\ref{eq:ineq_eO}) in Proposition~\ref{prop:ineq} 
that $S(g;C) \le 1$.
Hence, by using (\ref{eq:ineq_eC}), we conclude that $g \in \eCC(P)$.

Conversely, we show that $g \in \eCC(P)$ implies $\ePsi(g) \in \eOO(P)$.
Let $g \in \eCC(P)$ and put $f = \ePsi(g)$.
We need to prove that $T^+(f;C) \le 1$ for all $C \in \SC(P)$ with $\ctop C \in \max (P)$ 
and that $T^-(f;C) \le 1$ for all $C \in \MC(P)$.

Suppose $C = \{ v_1 \gtrdot v_2 \gtrdot \cdots \gtrdot v_r \} \in \SC(P)$ with $v_1 \in \max(P)$.
Then by definition
\begin{align*}
T^+(f;C)
 &=
2^{r-1} f(v_r) - \sum_{i=1}^{r-1} 2^{r-i-1} f(v_{r-i})
\\
 &=
2^{r-1} \left( g(v_r) + M(g;P_{< v_r}) \right)
 -
\sum_{i=1}^{r-1}
 2^{r-i-1} \left( g(v_{r-i}) + M(g;P_{< v_{r-i}}) \right).
\end{align*}
By using $x \le |x|$ and $-x \le |x|$, we see that
\[
T^+(f;C)
 \le
2^{r-1} \left( |g(v_r)| + M(g;P_{< v_r}) \right)
 +
\sum_{i=1}^{r-1}
 2^{r-i-1} \left( |g(v_{r-i})| - M(g;P_{< v_{r-i}}) \right).
\]
Then by using Lemma~\ref{lem:2-2} we obtain
\[
T^+(f;C)
 \le 
M(g;P_{\le v_1})
 =
\max \{ S(g;C') : C' \in \MC(P_{\le v_1}) \}.
\]
Since $S(g;C') \le 1$ for all $C' \in \MC(P_{\le v_1})$ by (\ref{eq:ineq_eC}), 
we have $T^+(f ; C) \le 1$.

Suppose $C = \{ v_1 \gtrdot v_2 \gtrdot \cdots \gtrdot v_r \} \in \MC(P)$. 
Then $v_1 \in \max(P)$ and $v_r \in \min(P)$.
It follows from the definition that
\begin{align*}
T^-(f;C)
 &=
- 2^{r-1} f(v_r) - \sum_{i=1}^{r-1} 2^{r-i-1} f(v_{r-i})
\\
 &=
- 2^{r-1} g(v_r)
 -
\sum_{i=1}^{r-1}
 2^{r-i-1} \left( g(v_{r-i}) + M(g;P_{< v_{r-i}}) \right).
\end{align*}
By using $x \le |x|$ and $-x \le |x|$, we see that
\[
T^-(f;C)
 \le
2^{r-1} |g(v_r)|
 +
\sum_{i=1}^{r-1}
 2^{r-i-1} \left( |g(v_{r-i})| - M(g;P_{< v_{r-i}}) \right).
\]
Since $\{ v_r \} \in \MC(P_{< v_{r-1}})$, 
we have $|g(v_r)| \le M(g;P_{< v_{r-1}})$.
Hence we have
\begin{align*}
&
T^-(f;C)
\\
&
 \le
2^{r-1} M(g;P_{< v_{r-1}})
 +
2^{r-2} \left( |g(v_{r-1})| - M(g;P_{< v_{r-1}}) \right)
\\
&\quad
 +
\sum_{i=2}^{r-1}
 2^{r-i-1} \left( |g(v_{r-i})| - M(g;P_{< v_{r-i}}) \right)
\\
&
=
2^{r-2} \left( |g(v_{r-1})| + M(g;P_{< v_{r-1}}) \right)
 +
\sum_{i=2}^{r-1}
 2^{r-i-1} \left( |g(v_{r-i})| - M(g;P_{< v_{r-i}}) \right).
\end{align*}
Now we can use Lemma~\ref{lem:2-2} and (\ref{eq:ineq_eO}) to obtain 
$T^-(f;C) \le M(g;P_{v_1}) \le 1$.

Therefore we conclude that $f \in \eOO(P)$.
\qed
\end{demo}

Here we show that the bijection $\ePhi : \eOO(P) \to \eCC(P)$ restricts to the bijection $\Phi : \OO(P) \to \CC(P)$.

\begin{prop}
\label{prop:e-transfer-O2C}
The restriction of the enriched transfer map $\ePhi$ to $\OO(P)$ 
coincides with the restriction of the transfer map $\Phi$ to $\OO(P)$.
\end{prop}

\begin{demo}{Proof}
Let $f \in \OO(P)$ and put $g = \Phi(f)$, $\tilde{g} = \ePhi(f)$.
By using the induction on the ordering of $P$, we prove
\begin{gather}
\label{eq:Phi1}
\max \{ f(w) : w \lessdot v \}
 =
\max \{ g(v_1) + \dots + g(v_r) : \{v_1 \gtrdot  \cdots \gtrdot v_r\} \in \MC(P_{<v}) \},
\\
\label{eq:Phi2}
\tilde{g}(v) = g(v) \ge 0.
\end{gather}
If $v$ is minimal in $P$, then $f(v) = g(v) = \tilde{g}(v)$.
If $v$ is not minimal in $P$ and $\{ w \in P : w \lessdot v \} = \{ w_1, \dots, w_k \}$, 
then it follows from the induction hypothesis for (\ref{eq:Phi2}) that
\begin{multline*}
\max \{ |\tilde{g}(v_1)| + \dots + |\tilde{g}(v_r)| : \{v_1 \gtrdot  \cdots \gtrdot v_r\} \in \MC(P_{< v}) \}
\\
= 
\max_{1 \le i \le k}
\left\{
 g(w_i) + \max \{ g(v_2) + \dots + v(v_r) : \{v_2 \gtrdot  \cdots \gtrdot v_r\} \in \MC(P_{< w_i}) \}
\right\}.
\end{multline*}
By using the induction hypothesis for (\ref{eq:Phi1}) and (\ref{eq:def_transfer}), we obtain
\begin{align*}
&
\max \{ |\tilde{g}(v_1)| + \dots + |\tilde{g}(v_r)| : \{v_1 \gtrdot  \cdots \gtrdot v_r\} \in \MC(P_{< v}) \}
\\
&\quad
 =
\max_{1 \le i \le k}
\big\{
 g(w_i) + \max \{ f(u_i) : u_i \lessdot w_i \}
\big\}
\\
&\quad
 =
\max_{1 \le i \le k} f(w_i)
 =
\max \{ f(w) : w \lessdot v \}.
\end{align*}
Hence, comparing (\ref{eq:e-transfer}) with (\ref{eq:def_transfer}), we obtain (\ref{eq:Phi1}) and (\ref{eq:Phi2}).
\qed
\end{demo}

\subsection{%
Left enriched $P$-partitions
}
\label{subsec:enrchiedbijection}
In this subsection, we use the enriched transfer map to find a bijection 
from left enriched $P$-partitions to lattice points of the dilated enriched order polytope.

Recall the definition of left enriched $P$-partition introduced by Petersen \cite{Petersen2007}.
A map $h : P \to \Int$ is called a \emph{left enriched $P$-partition} 
if it satisfies the following two conditions:
\begin{enumerate}
\item[(i)]
If $v \le w$, then $|h(v)| \le |h(w)|$;
\item[(ii)]
If $v \le w$ and $|h(v)| = |h(w)|$, then $h(w) \ge 0$.
\end{enumerate}
We denote by $\EE_m(P)$ the set of left enriched $P$-partitions $h : P \to \Int$ 
such that $|h(v)| \le m$ for all $v \in P$.
Note that $\eFF(P) = \EE_1(P)$.
Ohsugi--Tsuchiya \cite{OT1} gave an explicit bijection between $\EE_m(P)$ 
and $m \eCC(P) \cap \Int^P$.

\begin{prop}
\label{prop:bij_P2C}
(\cite[Theorem~0.2 and its proof]{OT1})
Let $\Pi : \EE_m(P) \to \Real^P$ be the map defined by
\begin{multline}
\label{eq:bij_P2C}
\left( \Pi(h) \right)(v)
\\
 =
\begin{cases}
 h(v) &\text{if $v$ is minimal in $P$,} \\
 h(v) - \max \{ |h(w)| : w \lessdot v \} &\text{if $v$ is not minimal in $P$ and $h(v) \ge 0$,} \\
 h(v) + \max \{ |h(w)| : w \lessdot v \} &\text{if $v$ is not minimal in $P$ and $h(v) < 0$.}
\end{cases}
\end{multline}
Then $\Pi$ gives a bijection from $\EE_m(P)$ to $m \eCC(P) \cap \Int^P$.
\end{prop}

By composing this bijection $\Pi$ with the inverse enriched transfer map 
$\ePsi : m \eCC(P) \cap \Int^P \to m \eOO(P) \cap \Int^P$, 
we obtain an explicit bijection from $\EE_m(P)$ to $m \eOO(P) \cap \Int^P$.

\begin{prop}
\label{prop:bij_P2O}
Let $\Theta : \EE_m(P) \to \Real^P$ be the map defined by
\begin{multline}
\label{eq:bij_P2O}
\left( \Theta(h) \right)(v)
\\
 =
\begin{cases}
 h(v) &\text{if $v$ is minimal in $P$ or $h(v) \ge 0$,} \\
 h(v) + 2 \max \{ |h(w)| : w \lessdot v \} &\text{if $v$ is not minimal in $P$ and $h(v) < 0$.}
\end{cases}
\end{multline}
Then $\Theta$ gives a bijection from $\EE_m(P)$ to $m \eOO(P) \cap \Int^P$.
\end{prop}

\begin{demo}{Proof}
We show that $\Theta = \ePsi \circ \Pi$.
Let $h \in \EE_m(P)$ and put $g = \Pi(h)$.
By comparing (\ref{eq:e-transfer}) with (\ref{eq:bij_P2C}) and (\ref{eq:bij_P2O}), 
it is enough to show
\begin{equation}
\label{eq:Pi1}
\max \{ S(g;C) : C \in \MC(P_{< v}) \}
 =
\max \{ |h(w)| : w \lessdot v \}.
\end{equation}
We proceed by induction on the ordering of $P$.
If $v$ is minimal in $P$, there is nothing to prove.
Suppose that $v$ is not minimal in $P$.
Since $h \in \EE_m(P)$, we have $|h(v)| \ge \max \{ |h(w)| : w \lessdot v \}$.
Then it follows from (\ref{eq:bij_P2C}) that 
\begin{equation}
\label{eq:Pi2}
|g(v)| = |h(v)| - \max \{ |h(w)| : w \lessdot v \}.
\end{equation}
If $\{ w \in P : w \lessdot v \} = \{ w_1, \dots, w_k \}$, then we have
\[
\max \{ S(g;C) : C \in \MC(P_{<v}) \}
 =
\max_{1 \le i \le k}
 \left\{ |g(w_i)| + \max \{ S(g;C') : C' \in \MC(P_{< w_i}) \} \right\}.
\]
By using (\ref{eq:Pi2}), we have
\begin{align*}
\max \{ S(g;C) : C \in \MC(P_{<v}) \}
 &=
\max_{1 \le i \le k}
 \left \{ |g(w_i)| + \max \{ |h(u_i)| : u_i \lessdot w_i \} \right\}
\\
 &=
\max_{1 \le i \le k} \{ |h(w_i)| \}, 
\end{align*}
from which (\ref{eq:Pi1}) follows.
\qed
\end{demo}

\subsection{%
Vertices of enriched poset polytopes
}
\label{sec:vertices}

In this subsection, we determine the vertex sets of the enriched order polytope $\eOO(P)$ 
and the enriched chain polytope $\eCC(P)$.

In order to state the result, we need a partial ordering $\preceq$ on $\eFF(P)$ or $\eAC(P)$.
For $f$, $f' \in \eFF(P)$ (or $\eAC(P)$), 
we write $f \preceq f'$ if $\supp(f) \subset \supp(f')$ and $f|_{\supp(f)} = f'|_{\supp(f)}$.

\begin{example}
\label{ex:vertex}
If $\Lambda = \{ u,v,w \}$ is the three-element chain with covering relations 
$u \lessdot w$ and $v \lessdot w$, 
then the Hasse diagrams of $\eFF(P)$ and $\eAC(P)$ with respect to $\preceq$ are shown 
in Figures~\ref{fig:Hasse2} and \ref{fig:Hasse3} respectively. %
\begin{figure}[htb]
\centering
					\begin{tikzpicture}[scale=1]
\coordinate (A) at (3,0); 
\coordinate (B) at (0,1.5);
\coordinate (C) at (6,1.5);
\coordinate (D) at (-4,3);
\coordinate (E) at (-1.5,3);
\coordinate (F) at (1,3);
\coordinate (G) at (3.5,3);
\coordinate (H) at (-4,4.5);
\coordinate (I) at (-1.5,4.5);
\coordinate (J) at (1,4.5);
\coordinate (K) at (3.5,4.5);

\coordinate (A2) at (3,0.3); 
\coordinate (B1) at (0,1.2);
\coordinate (B2) at (0,1.8);
\coordinate (C1) at (6,1.2);
\coordinate (C2) at (6,1.8);
\coordinate (D1) at (-4,2.7);
\coordinate (D2) at (-4,3.3);
\coordinate (E1) at (-1.5,2.7);
\coordinate (E2) at (-1.5,3.3);
\coordinate (F1) at (1,2.7);
\coordinate (F2) at (1,3.3);
\coordinate (G1) at (3.5,2.7);
\coordinate (G2) at (3.5,3.3);
\coordinate (H1) at (-4,4.2);
\coordinate (I1) at (-1.5,4.2);
\coordinate (J1) at (1,4.2);
\coordinate (K1) at (3.5,4.2);

\node at (A) {$(0,0,0)$};
\node at (B)  {$(0,0,1)$};
\node at (C)  {$(0,0,-1)$};
\node at (D)  {$(1,0,1)$};
\node at (E)  {$(0,-1,1)$};
\node at (F)  {$(-1,0,1)$};
\node at (G)  {$(0,1,1)$};
\node at (H)  {$(1,1,1)$};
\node at (I)  {$(1,-1,1)$};
\node at (J)  {$(-1,-1,1)$};
\node at (K)  {$(-1,1,1)$};

\draw (A2)--(B1);
\draw (A2)--(C1);
\draw (B2)--(D1);
\draw (B2)--(E1);
\draw (B2)--(F1);
\draw (B2)--(G1);
\draw (D2)--(H1);
\draw (D2)--(I1);
\draw (E2)--(I1);
\draw (E2)--(J1);
\draw (F2)--(J1);
\draw (F2)--(K1);
\draw (G2)--(K1);
\draw (G2)--(H1);
\end{tikzpicture}
\caption{Hasse diagram of $(\eFF(\Lambda),\preceq)$}
\label{fig:Hasse2}
\end{figure}
\begin{figure}[htb]
\centering
					\begin{tikzpicture}[scale=1]
\coordinate (A) at (0,0); 
\coordinate (B) at (-6,1.5);
\coordinate (C) at (-3.6,1.5);
\coordinate (D) at (-1.2,1.5);
\coordinate (E) at (1.2,1.5);
\coordinate (F) at (3.6,1.5);
\coordinate (G) at (6.0,1.5);
\coordinate (H) at (-6,3);
\coordinate (I) at (-3.6,3);
\coordinate (J) at (-1.2,3);
\coordinate (K) at (1.2,3);

\coordinate (A2) at (0,0.3); 
\coordinate (B1) at (-6,1.2);
\coordinate (C1) at (-3.6,1.2);
\coordinate (D1) at (-1.2,1.2);
\coordinate (E1) at (1.2,1.2);
\coordinate (B2) at (-6,1.8);
\coordinate (C2) at (-3.6,1.8);
\coordinate (D2) at (-1.2,1.8);
\coordinate (E2) at (1.2,1.8);
\coordinate (F1) at (3.6,1.2);
\coordinate (G1) at (6.0,1.2);
\coordinate (H1) at (-6,2.7);
\coordinate (I1) at (-3.6,2.7);
\coordinate (J1) at (-1.2,2.7);
\coordinate (K1) at (1.2,2.7);

\node at (A) {$(0,0,0)$};
\node at (B)  {$(1,0,0)$};
\node at (C)  {$(0,-1,0)$};
\node at (D)  {$(-1,0,0)$};
\node at (E)  {$(0,1,0)$};
\node at (F)  {$(0,0,1)$};
\node at (G)  {$(0,0,-1)$};
\node at (H)  {$(1,1,0)$};
\node at (I)  {$(1,-1,0)$};
\node at (J)  {$(-1,-1,0)$};
\node at (K)  {$(-1,1,0)$};

\draw (A2)--(B1);
\draw (A2)--(C1);
\draw (A2)--(D1);
\draw (A2)--(E1);
\draw (A2)--(F1);
\draw (A2)--(G1);
\draw (B2)--(H1);
\draw (C2)--(I1);
\draw (D2)--(J1);
\draw (E2)--(K1);
\draw (B2)--(I1);
\draw (C2)--(J1);
\draw (D2)--(K1);
\draw (E2)--(H1);
\end{tikzpicture}
\caption{Hasse diagram of $(\eAC(\Lambda),\preceq)$}
\label{fig:Hasse3}
\end{figure}
The enriched order polytope $\eOO(\Lambda)$ is the pyramid with five vertices 
$(1,1,1), (1,-1,1), (-1,-1,1), (-1,1,1)$ and $(0,0,-1)$, 
while the enriched chain polytope $\eCC(\Lambda)$ is the bipyramid with six vertices 
$(1,1,1), (1,-1,1), (-1,-1,1), (-1,1,1), (0,0,1)$ and $(0,0,-1)$.
\end{example}

\begin{prop}
\label{prop:e-vertex}
\begin{enumerate}
\item[(a)]
A point $f \in \eFF(P)$ is a vertex of $\eOO(P)$ 
if and only if $f$ is maximal with respect to the ordering $\preceq$.
\item[(b)]
A point $f \in \eAC(P)$ is a vertex of $\eCC(P)$ 
if and only if $f$ is maximal with respect to the ordering $\preceq$.
\end{enumerate}
\end{prop}

Note that $f \in \eAC(P)$ is maximal with respect to $\preceq$ 
if and only if $\supp(f)$ is a maximal antichain.

\begin{demo}{Proof}
(a)
Let $f$ be a maximal element of $\eFF_P$ with respect to $\preceq$.
Assume to the contrary that $f$ is not a vertex of $\eOO(P)$.
Then there exist elements $g_1, \dots, g_r \in \eFF(P)$ 
and positive real numbers $\lambda_1, \dots, \lambda_r$ such that $g_i \neq f$ 
and
\[
f = \sum_{i=1}^r \lambda_i g_i,
\quad
\sum_{i=1}^r \lambda_i = 1.
\]
Considering the value at $v \in P$, we have
\[
\sum_{i=1}^r \lambda_i g_i(v) = f(v) = \sum_{i=1}^r \lambda_i f(v).
\]
If $f(v) = 1$, then we see that
$\sum_{i=1}^r \lambda_i \left( 1 - g_i(v) \right) = 0$.
Since $\lambda_i > 0$ and $1 - g_i(v) \ge 0$, 
we obtain $g_i(v) = 1$ for all $i$.
By a similar reasoning, we see that, if $f(v) = -1$, then we have $g_i(v) = -1$ for all $i$.
Hence we have $\supp(f) \subset \supp(g_i)$ and $f|_{\supp(f)} = g_i|_{\supp(f)}$.
Since $f$ is maximal with respect to $\preceq$, 
we have $f = g_i$, which contradicts to the assumption $g_i \neq f$.
Therefore $f$ is a vertex of $\eOO(P)$.

Conversely, suppose that $f$ is not maximal with respect to $\preceq$.
Then there exists $g \in \eFF(P)$ such that 
$\supp(f) \subsetneq \supp(g)$ and $f|_{\supp(f)} = g|_{\supp(f)}$.
We take a maximal element $u$ of $\supp(g) \setminus \supp(f)$ and define $f'$, $f'' : P \to \Real$ by
$$
f'(v)
 = 
\begin{cases}
 f(v) = g(v) &\text{if $v \in \supp(f)$,} \\
 g(u) &\text{if $v = u$,} \\
 0 &\text{otherwise,}
\end{cases}
\quad
f''(v)
 = 
\begin{cases}
 f(v) = g(v) &\text{if $v \in \supp(f)$,} \\
 -g(u) &\text{if $v = u$,} \\
 0 &\text{otherwise}.
\end{cases}
$$
Then $\supp(f') = \supp(f'') = \supp(f) \sqcup \{ u \}$ is an order filter of $P$ 
and $u$ is an minimal element of $\supp(f') = \supp(f'')$.
Hence $f' \in \eFF(P)$.
Since $f = (f' + f'')/2$, we see that $f$ is not a vertex of $\eOO(P)$.

(b)
Similar to (a).
\qed
\end{demo}

\begin{remark}
\begin{enumerate}
\item[(1)]
A characterization of the vertex set of the enriched chain polytope $\eCC(P)$ 
is also given in \cite[Section~7]{KOS}.
\item[(2)]
In general, the image $\ePhi(f)$ of a vertex $f$ of $\eOO(P)$ under the enriched transfer map $\ePhi$ 
is not a vertex of $\eCC(P)$, 
and the number of vertices of $\eOO(P)$ is different from that of $\eCC(P)$
(see Example~\ref{ex:vertex}).
\end{enumerate}
\end{remark}

\section{%
Triangulations
}
\label{sec:triangulation}

In this section we prove Theorem~\ref{thm:e-triangulation}, 
which describes triangulations of enriched order and chain polytopes.

\subsection{%
Poset structure on $\eFF(P)$
}

We introduce a partial ordering $\ge$ on $\eFF(P)$, 
which is an extension of the inclusion ordering on the set of order filters of $P$.
Note that this ordering $\ge$ is different from the ordering $\succeq$ used in 
Section~\ref{sec:vertices}.

\begin{definition}
\label{def:ordering}
For $f$, $g \in \eFF(P)$, we write $f > g$ if the following three conditions hold:
\begin{enumerate}
\item[(i)]
$\supp (f) \supsetneq \supp (g)$;
\item[(ii)]
$f(v) \ge g(v)$ for any $v \in \supp(g)$;
\item[(iii)]
If $v \in \supp(g)$ and $v$ is minimal in $\supp(f)$, then $f(v) = g(v)$.
\end{enumerate}
Also we write $f \ge g$ if $f=g$ or $f>g$.
\end{definition}

The following lemma is obvious, but will be used in several places.

\begin{lemma}
\label{lem:minimal}
If $F \supset G$ are order filters of $P$ and $v \in G$ is minimal in $F$, 
then $v$ is minimal in $G$.
\end{lemma}

By using this lemma, we can prove that $\eFF(P)$ is equipped with a poset structure 
with respect to the binary relation $\ge$.

\begin{lemma}
\label{lem:ordering}
The binary relation $\ge$ given in Definition~\ref{def:ordering} is a partial ordering 
on $\eFF(P)$.
\end{lemma}

\begin{demo}{Proof}
It is enough to show the transitivity.
Let $f$, $g$, $h \in \eFF(P)$ satisfy $ f > g$ and $g > h$.
Then it is clear that $\supp(f) \supsetneq \supp(h)$ and $f(v) \ge h(v)$ for any $v \in \supp(h)$.
Since $\supp(f) \supset \supp(g) \supset \supp(h)$, 
it follows from Lemma~\ref{lem:minimal} that, if $v \in \supp(h)$ is minimal in $\supp(f)$, 
then we have $f(v) = g(v) = h(v)$.
\qed
\end{demo}

\begin{example}
Let $\Lambda$ be the three-element poset on $\{ u, v, w \}$ 
with covering relations $u \lessdot w$ and $v \lessdot w$.
Figure~\ref{fig:Hasse1} shows the Hasse diagram of 
$(\eFF(\Lambda), \ge)$.
\begin{figure}
\centering
					\begin{tikzpicture}[scale=1]
\coordinate (A) at (0,0); 
\coordinate (B) at (-1.5,1.5);
\coordinate (C) at (1.5,1.5);
\coordinate (D) at (-4.5,3);
\coordinate (E) at (-1.5,3);
\coordinate (F) at (1.5,3);
\coordinate (G) at (4.5,3);
\coordinate (H) at (-4.5,4.5);
\coordinate (I) at (-1.5,4.5);
\coordinate (J) at (1.5,4.5);
\coordinate (K) at (4.5,4.5);

\coordinate (A2) at (0,0.3); 
\coordinate (B1) at (-1.5,1.2);
\coordinate (C1) at (1.5,1.2);
\coordinate (D1) at (-4.5,2.7);
\coordinate (E1) at (-1.5,2.7);
\coordinate (F1) at (1.5,2.7);
\coordinate (G1) at (4.5,2.7);
\coordinate (H1) at (-4.5,4.2);
\coordinate (I1) at (-1.5,4.2);
\coordinate (J1) at (1.5,4.2);
\coordinate (K1) at (4.5,4.2);
\coordinate (B2) at (-1.5,1.8);
\coordinate (C2) at (1.5,1.8);
\coordinate (D2) at (-4.5,3.3);
\coordinate (E2) at (-1.5,3.3);
\coordinate (F2) at (1.5,3.3);
\coordinate (G2) at (4.5,3.3);

\node at (A) {$(0,0,0)$};
\node at (B)  {$(0,0,1)$};
\node at (C)  {$(0,0,-1)$};
\node at (D)  {$(1,0,1)$};
\node at (E)  {$(0,-1,1)$};
\node at (F)  {$(-1,0,1)$};
\node at (G)  {$(0,1,1)$};
\node at (H)  {$(1,1,1)$};
\node at (I)  {$(1,-1,1)$};
\node at (J)  {$(-1,-1,1)$};
\node at (K)  {$(-1,1,1)$};

\draw (A2)--(B1);
\draw (A2)--(C1);
\draw (B2)--(D1);
\draw (B2)--(E1);
\draw (B2)--(F1);
\draw (B2)--(G1);
\draw (C2)--(D1);
\draw (C2)--(E1);
\draw (C2)--(F1);
\draw (C2)--(G1);
\draw (D2)--(H1);
\draw (D2)--(I1);
\draw (E2)--(I1);
\draw (E2)--(J1);
\draw (F2)--(J1);
\draw (F2)--(K1);
\draw (G2)--(K1);
\draw (G2)--(H1);
\end{tikzpicture}
\caption{Hasse diagram of $(\eFF(\Lambda), \ge)$}
\label{fig:Hasse1}
\end{figure}
\end{example}

We collect several properties of this partial ordering on $\eFF(P)$.

\begin{prop}
\label{prop:ordering}
The resulting poset $\eFF(P)$ has the following properties.
\begin{enumerate}
\item[(a)]
For order filters $F$ and $G$, we have $F \supset G$ if and only if $\chi_F \ge \chi_G$ 
in $\eFF(P)$, where $\chi_S$ is the characteristic function of $S$.
\item[(b)]
The zero map $0$ is the unique minimal element of $\eFF(P)$.
\item[(c)]
If $f$ covers $g$ in $\eFF(P)$, then $\# \supp(f) = \# \supp(g) + 1$.
\item[(d)]
If $f$ is a maximal element in $\eFF(P)$, then $\supp(f) = P$.
\item[(e)]
All maximal chains of $\eFF(P)$ have the same length $d = \# P$.
\end{enumerate}
\end{prop}

\begin{demo}{Proof}
(a) and (b) are obvious.

(c)
It is enough to show that, 
if $f>g$, then there exists $h \in \eFF(P)$ such that $f \ge h  > g$ 
and $\# \supp(h) = \# \supp(g) + 1$.

Since $\supp(f) \supsetneq \supp(g)$ and they are order filters of $P$, 
there exists $u \in \supp(f)$ such that $\supp(g) \cup \{ u \}$ is an order filter of $P$.
Then we define $h : P \to \{ 1, 0, -1 \}$ by putting
\[
h(v)
 =
\begin{cases}
 f(v) &\text{if $v \in \supp(g) \cup \{ u \}$,} \\
 0 &\text{otherwise.}
\end{cases}
\]
We see that $h \in \eFF(P)$, $\supp(h) = \supp(g) \cup \{ u \}$, 
and $f \ge h > g$.

(d)
Suppose that $\supp(g) \neq P$.
Since $\supp(g)$ is a proper order filter of $P$, there exists $u \not\in \supp(g)$ 
such that $\supp(g) \cup \{ u \}$ is an order filter.
Define $f : P \to \{ 1, 0, -1 \}$ by putting
\[
f(v)
 =
\begin{cases}
 1 &\text{if $v=u$ or $v$ covers $u$,} \\
 g(v) &\text{otherwise.}
\end{cases}
\]
Then we have $f \in \eFF(P)$, $\supp(f) = \supp(g) \cup \{ u \}$ and $f > g$.

(e) follows from (b), (c) and (d).
\qed
\end{demo}

Next we consider chains in the poset $\eFF(P)$.

\begin{definition}
\label{def:supp-sgn}
Given a chain $K = \{ f_1 > f_2 > \dots > f_k \}$ of $\eFF(P)$, 
we define its \emph{support} $\supp(K)$ and \emph{signature} $\sgn(K)$ as follows.
The support $\supp(K)$ is the chain $\{ \supp(f_1) \supsetneq \supp(f_2) \supsetneq \dots \supsetneq \supp(f_k) \}$ 
of order filters.
The signature $\sgn(K)$ is the map $\varphi : P \to \{ 1, 0, -1 \}$ given by
\begin{enumerate}
\item[(i)]
If $v$ is not minimal in $\supp(f_i)$ for any $i$, 
then $\varphi(v) = 0$;
\item[(ii)]
If $v$ is minimal in $\supp(f_i)$ for some $i$, 
then $\varphi(v) = f_i(v)$.
\end{enumerate}
\end{definition}

The following lemma guarantees that the definition of $\varphi(v)$ in the case (ii) 
is independent of the choice of $i$.

\begin{lemma}
\label{lem:sgn}
Let $K = \{ f_1 > f_2 > \dots > f_k \}$ be a chain of $\eFF(P)$.
If $v$ is minimal in both $\supp(f_i)$ and $\supp(f_j)$, 
then we have $f_i(v) = f_j(v)$.
\end{lemma}

\begin{demo}{Proof}
We may assume $i<j$. Then $f_i > f_j$ and $\supp(f_i) \supset \supp(f_j)$.
Since $v \in \supp(f_j)$ and minimal in $\supp(f_i)$, we have $f_i(v) = f_j(v)$ 
by the condition (iii) in Definition~\ref{def:ordering}.
\qed
\end{demo}

A key property of support and signature is the following.

\begin{prop}
\label{prop:bijection}
Let $X(P)$ be the set of all chains of $\eFF(P)$ (including the empty chain), 
and $Y(P)$ the set of all pairs $(C, \varphi)$ of chains 
$C = \{ F_1 \supsetneq F_2 \supsetneq \dots \supsetneq F_k \}$ of order filters of $P$ 
and maps $\varphi : P \to \{ 1, 0, -1 \}$ satisfying
\begin{equation}
\label{eq:bijection}
\supp(\varphi) = \bigcup_{i=1}^k \min F_i,
\end{equation}
where $\min F_i$ is the set of minimal elements of $F_i$.
Then the map $X(P) \ni K \mapsto (\supp(K), \sgn(K)) \in Y(P)$ is a bijection.
In particular, maximal chains in $\eFF(P)$ are in bijection with pairs $(C,\varphi)$ 
of maximal chains $C$ of order filters and maps $\varphi : P \to \{ 1, -1 \}$.
\end{prop}

It follows that the number of maximal chains in $\eFF(P)$ is equal to 
$2^d e(P)$, where $d = \# P$ and $e(P)$ is the number of linear extensions of $P$.

\begin{demo}{Proof}
It follows from Definition~\ref{def:supp-sgn} that $(\supp(K), \sgn(K)) \in Y(P)$ for $K \in X(P)$.

Given a chain $C = \{ F_1 \supsetneq \dots \supsetneq F_k \}$ of order filters 
and a map $\varphi : P \to \{ 1, 0, -1 \}$ satisfying (\ref{eq:bijection}), 
we define $f_1, \cdots, f_k \in \Real^P$ by
\[
f_i(v)
 = 
\begin{cases}
 1 &\text{if $v \in F_i$ and $v$ is not minimal in $F_i$,} \\
 \varphi(v) &\text{if $v \in F_i$ and $v$ is minimal in $F_i$,} \\
 0 &\text{if $v \not\in F_i$.}
\end{cases}
\]
Then we see that $f_i \in \eFF(P)$ and $\supp(f_i) = F_i$.

We show that $f_i > f_{i+1}$ for $1 \le i \le k-1$.
Firstly one has $\supp(f_i) = F_i \supsetneq F_{i+1} = \supp(f_{i+1})$.
Secondly we check that $f_i(v) \ge f_{i+1}(v)$ for $v \in \supp(f_{i+1})$.
Since $v \in \supp(f_{i+1}) \subset \supp(f_i)$, we have $f_i(v)$, $f_{i+1}(v) \in \{ 1, -1 \}$, 
and there is nothing to prove in the case $f_i(v) = 1$.
If $f_i(v) = -1$, then $v$ is minimal in $\supp(f_i)$, so $v$ is minimal in $\supp(f_{i+1})$ 
by Lemma~\ref{lem:minimal}.
Then we have $\varphi(v) = -1$ and $f_{i+1}(v) = -1 = f_i(v)$.
Lastly, if $v \in \supp(f_{i+1})$ and $v$ is minimal in $\supp(f_i)$, 
then $v$ is minimal in $\supp(f_{i+1})$ by Lemma~\ref{lem:minimal} 
and $f_i(v) = \varphi(v) = f_{i+1}(v)$. 

Therefore $K = \{ f_1 > f_2 > \dots > f_k \}$ is a chain in $\eFF(P)$, 
and $\supp(K) = \{ F_1 \supsetneq F_2 \supsetneq \dots \supsetneq F_k \}$, 
$\sgn(K) = \varphi$.
\qed
\end{demo}

\subsection{%
Triangulation of $\eCC(P)$
}

In this subsection, we use the triangulation of $\CC(P)$ given in Theorem~\ref{thm:triangulation} 
to construct a unimodular triangulation of $\eCC(P)$.
We transfer this triangulation of $\eCC(P)$ to $\eOO(P)$ 
via the inverse enriched transfer map $\ePsi$ in the next subsection.

A (lattice) \emph{triangulation} of a lattice polytope $\mathcal{P} \subset \Real^d$ 
of dimension $d$ is a finite collection $\Delta$ of (lattice) simplices such that
\begin{enumerate}
\item[(i)]
every face of a member of $\Delta$ is in $\Delta$,
\item[(ii)]
the union of the simplices in $\Delta$ is $\mathcal{P}$, and
\item[(iii)]
any two elements of $\Delta$ intersect in a common (possibly empty) face.
\end{enumerate}
We say that a triangulation $\Delta$ is \emph{unimodular} 
if all maximal faces of $\Delta$ are unimodular, i.e., have the Euclidean volume $1/d!$.

Recall that the simplices $S_C=\conv \{ \chi_{F} : F \in C\}$ and $T_C=\conv\{ \Phi (\chi_{F}): F \in C\}$ of the triangulation given in Theorem~\ref{thm:triangulation} are 
described as follows.

\begin{prop}
(Stanley \cite[Section~5]{Stanley1986})
\label{prop:simplex}
If $C = \{ F_1 \supsetneq F_2 \supsetneq \dots \supsetneq F_k \}$ is a chain of order filters of $P$, then we have
\begin{equation}
\label{eq:simplex1}
S_C
 = 
\left\{
 \begin{array}{l}
 f \in \Real^P :
 \\
 \quad
 \text{(i) $f$ is constant on the subsets $P \setminus F_1, F_1 \setminus F_2, \dots, F_{k-1} \setminus F_k, F_k$},
 \\
 \quad
 \text{(ii) $0 = f(P \setminus F_1) \le f(F_1 \setminus F_2) \le \cdots \le f(F_{k-1} \setminus F_k) \le f(F_k) = 1$}.
 \end{array}
\right\},
\end{equation}
and
\begin{equation}
\label{eq:simplex2}
T_C = \Phi(S_C).
\end{equation}
\end{prop}

If $C = \{ F_1 \supsetneq F_2 \supsetneq \dots \supsetneq F_k \}$ 
is a chain of order filters of $P$, then 
$\chi_C = \{ \chi_{F_1} > \chi_{F_2} > \dots > \chi_{F_k} \}$ is a chain in $\eFF(P)$ 
by Proposition~\ref{prop:ordering} (a), and
\[
S^{(e)}_{\chi_C} = S_C,
\quad
T^{(e)}_{\chi_C} = T_C.
\]
First we show that any $T^{(e)}_K =\conv \ePhi(K)$ is obtained from $T_C$ 
by a composition of reflections.
For $\varphi : P \to \{ 1, 0, -1 \}$, we define a linear map $R_\varphi : \Real^P \to \Real^P$ by 
\[
\left( R_\varphi g \right)(v)
 = 
\begin{cases}
 g(v) &\text{if $\varphi(v) = 1$ or $0$,} \\
 -g(v) &\text{if $\varphi(v) = -1$.}
\end{cases}
\]
The linear map $R_\varphi$ is a composition of reflections along coordinate hyperplanes.

\begin{prop}
\label{prop:simplex-reflection}
For a chain $K$ in $\eFF_P$, we obtain
\begin{equation}
\label{eq:simplex-reflection}
T^{(e)}_K = R_{\sgn(K)} (T_{\supp(K)}).
\end{equation}
\end{prop}

\begin{demo}{Proof}
Let $K = \{ f_1 > \dots > f_k \}$ and 
put $C = \supp(K) = \{ F_1 \supsetneq \dots \supsetneq F_k \}$ ($F_i = \supp(f_i)$) 
and $\varphi = \sgn(K)$.
Since $T_C = \conv \Phi(\chi_C)$, we have
\[
R_\varphi T_C
 =
R_\varphi ( \conv \ePhi(\chi_C) )
 =
\conv ( R_\varphi (\ePhi (\chi_C)) ).
\]
Hence it is enough to show that $R_\varphi (\ePhi (\chi_{F_i}) ) = \ePhi(f_i)$ for each $i$.

By the definition of the enriched transfer map, we have
\begin{align*}
\ePhi(\chi_{F_i})(v)
 &=
\begin{cases}
 1 &\text{if $v$ is minimal in $F_i$,} \\
 0 &\text{otherwise,}
\end{cases}
\\
\ePhi(f_i)(v)
 &=
\begin{cases}
 f_i(v) &\text{if $v$ is minimal in $F_i$,} \\
 0 &\text{otherwise.}
\end{cases}
\end{align*}
On the other hand, it follows from the definition of $\varphi = \sgn(K)$ that
\[
\varphi(v)
 =
\begin{cases}
 f_i(v) &\text{if $v$ is minimal in some $\supp(f_i)$,} \\
 0 &\text{otherwise.}
\end{cases}
\]
Hence we obtain $R_\varphi (\ePhi (\chi_{F_i}) ) = \ePhi(f_i)$.
\qed
\end{demo}

In order to prove Theorem~\ref{thm:e-triangulation} (b), we prepare several lemmas.
Given $\varphi \in \{ 1, 0, -1 \}^P$, we put
\[
V_\varphi
 =
\left\{
 g \in \Real^P :
\begin{array}{l}
\text{(i) if $\varphi(v) = 1$, then $g(v) \ge 0$,} \\
\text{(ii) if $\varphi(v) = 0$, then $g(v) = 0$} \\
\text{(iii) if $\varphi(v) = -1$, then $g(v) \le 0$}
\end{array}
\right\}.
\]
For $\ep \in \{ 1, -1 \}^P$, we put
\[
\eCC_\ep(P) = \eCC(P) \cap V_\ep,
\quad
\eAC_\ep(P) = \eAC(P) \cap V_\ep.
\]
Since $\Real^P = \bigcup_{\ep \in \{ 1, -1  \}^P} V_\ep$, we have
\[
\eCC(P) = \bigcup_{\ep \in \{ 1, -1 \}^P} \eCC_\ep(P),
\quad
\eAC(P) = \bigcup_{\ep \in \{ 1, -1 \}^P} \eAC_\ep(P).
\]

\begin{lemma}
\label{lem:3-1}
(\cite[lemma~1.1]{OT1})
For $\ep \in \{ 1, -1 \}^P$, we have
\[
\eCC_\ep(P) = \conv ( \eAC_\ep(P) ) = R_\ep ( \CC(P) ).
\]
\end{lemma}

\begin{demo}{Proof}
The first equality is proved in \cite[Lemma~1.1]{OT1}.
We prove the second equality.
Let $\ep_0$ be the map given by $\ep_0(v) = 1$ for all $v \in P$.
Then $\eAC_{\ep_0}(P) = \AC(P)$ and $\eCC_{\ep_0}(P) = \conv ( \AC(P) ) = \CC(P)$.
Since $\eAC_\ep(P) = R_\ep ( \eAC_{\ep_0}(P) ) = R_\ep ( \AC(P) )$, we have
\[
\eCC_\ep(P)
 =
\conv ( R_\ep ( \AC(P) ) )
 =
R_\ep ( \conv ( \AC(P) ) )
 =
R_\ep ( \CC(P) ).
\]
\qed
\end{demo}

\begin{lemma}
\label{lem:3-2}
Suppose that $\varphi \in \{ 1, 0, -1 \}^P$ and $\ep \in \{ 1, -1 \}^P$ satisfy 
$\varphi|_{\supp(\varphi)} = \ep|_{\supp(\varphi)}$.
Then we have
\begin{enumerate}
\item[(a)]
$V_\varphi \subset V_\ep$.
\item[(b)]
$R_\varphi|_{V_{|\varphi|}} = R_\ep|_{V_{|\varphi|}}$, 
where $|\varphi|$ is defined by $|\varphi|(v) = |\varphi(v)|$.
\end{enumerate}
\end{lemma}

\begin{demo}{Proof}
(a)
Let $g \in V_\varphi$.
If $\ep(v) = 1$, then $\varphi(v) = 1$ or $0$ and $g(v) \ge 0$.
If $\ep(v) = -1$, then $\varphi(v) = -1$ or $0$ and $g(v) \le 0$.
Hence $g \in V_\ep$.

(b)
Let $g \in V_{|\varphi|}$.
If $v \in \supp(\varphi)$, then we have $\varphi(v) = \ep(v)$ and
$(R_\varphi g)(v) = \varphi(v) g(v) = \ep(v) g(v) = (R_\ep g)(v)$.
If $v \not\in \supp(\varphi)$, then we have $\varphi(v) = g(v) = 0$, thus
$(R_\varphi g)(v) = g(v) = 0$ and $(R_\ep g)(v) = \ep(v) g(v) = 0$.
\qed
\end{demo}

\begin{lemma}
\label{lem:3-3}
\begin{enumerate}
\item[(a)]
Let $C = \{ F_1 \supsetneq \dots \supsetneq F_k \}$ be a chain 
of order filters of $P$.
If $\varphi \in \{ 1, 0, -1 \}^P$ satisfies $\supp(\varphi) = \bigcup_{i=1}^k \min F_i$, 
then $T_C \subset V_{|\varphi|}$.
\item[(b)]
If $K$ is a chain in $\eFF(P)$, then we have $T^{(e)}_K \subset V_{\sgn(K)}$.
\item[(c)]
If $K$ is a chain in $\eFF(P)$ 
and $\ep \in \{ 1, -1 \}^P$ satisfies $\sgn(K)|_{\supp (\sgn(K))} = \ep|_{\supp (\sgn(K))}$, 
then we have $T^{(e)}_K = R_\ep T_{\supp(K)}$.
\end{enumerate}
\end{lemma}

\begin{demo}{Proof}
(a)
Let $g \in T_C$.
It is enough to show that $\varphi(v) = 0$ implies $g(v) = 0$.
By Theorem~\ref{thm:triangulation} (a), 
there exists $f \in S_C$ such that $g = \Phi(f)$.
Let $i$ be the largest index such that $v \in F_i$, where we use the convention $F_0 = P$.
If $i=0$, then $f(v) = 0$ and $g(v) = 0$.
Suppose that $i \ge 1$ and $\varphi(v) = 0$.
Then $v \in F_i \setminus F_{i+1}$ by the maximality of $i$.
Since $v$ is not minimal in $F_i$, there exists $w \in F_i$ such that $w \lessdot v$.
If $w \in F_{i+1}$, then $v \in F_{i+1}$ (since $F_{i+1}$ is an order filter), 
which contradicts to the maximality of $i$.
Hence we have $w \in F_i \setminus F_{i+1}$.
Then by (\ref{eq:simplex1}), we have $f(v) = f(w)$.
Therefore $g(v) = (\Phi(f))(v) = f(v) - \max \{ f(u) : u \lessdot v \} = f(v) - f(w) = 0$.

(b)
Let $C = \supp(K)$ and $\varphi = \sgn(K)$.
By (a), we have $T_C \subset V_{|\varphi|}$. 
Since $R_\varphi V_{|\varphi|} = V_\varphi$, 
we obtain $T_K = R_\varphi (T_C) \subset V_\varphi$.

(c) follows from Proposition~\ref{prop:simplex-reflection}, (a) and Lemma~\ref{lem:3-2} (b).
\qed
\end{demo}

Note that $\ePhi$ gives a bijection between $\eFF(P)$ and $\eAC(P)$ (Proposition~\ref{prop:e-transfer-F2A} (c)), 
and that $R_\varphi$ preserves $\eAC(P)$ for any $\varphi \in \{ 1, 0, -1 \}^P$.

\begin{lemma}
\label{lem:3-4}
Given $f_1$, $f_2 \in \eFF(P)$ and $\varphi \in \{ 1, 0, -1 \}^P$, 
we define $f'_1$, $f'_2 \in \eFF(P)$ by the condition
\[
R_\varphi \ePhi(f_1) = \ePhi(f'_1),
\quad
R_\varphi \ePhi(f_2) = \ePhi(f'_2).
\]
Then $f_1 > f_2$ implies $f'_1 > f'_2$.
\end{lemma}

\begin{demo}{Proof}
We may assume that there exists a unique $u \in P$ such that $\varphi(u) = -1$, 
i.e., 
\[
(R_\varphi g)(v)
 = 
\begin{cases}
 g(v) &\text{if $v \neq u$,} \\
 -g(v) &\text{if $v=u$.}
\end{cases}
\]
Then it follows from Proposition~\ref{prop:e-transfer-F2A} that, if $R_\varphi (\ePhi(f)) = \ePhi(f')$, then
\[
f'(v)
 =
\begin{cases}
 -f(v) &\text{if $u \in \min (\supp(f))$ and $v=u$,} \\
 f(v) &\text{otherwise,}
\end{cases}
\]
and $\supp(f') = \supp(f)$.

Now we assume that $f_1 > f_2$.
Then it is enough to prove the following two claims:
\begin{enumerate}
\item[(1)]
If $u \in \supp(f'_2)$, then $f'_1(u) \ge f'_2(u)$.
\item[(2)]
If $u \in \supp(f'_2)$ and $u$ is minimal in $\supp(f'_1)$, then $f'_1(u) = f'_2(u)$.
\end{enumerate}

First we prove (1) by dividing into four cases.
If $u \in \min (\supp(f_1))$ and $u \in \min (\supp(f_2))$, 
then we have $f_1(u) = f_2(u)$, thus $f'_1(u) = - f_1(u) = - f_2(u) = f'_2(u)$.
If $u \in \min (\supp(f_1))$ and $u \not\in \min (\supp(f_2))$, 
then it follows from Lemma~\ref{lem:minimal} that $u \not\in \supp(f_2)$, 
which contradicts to the assumption $u \in \supp(f'_2) = \supp(f_2)$.
If $u \not\in \min (\supp(f_1))$ and $u \in \min (\supp(f_2))$, 
then we have $f'_1(u) = f_1(u) = 1$, thus $f'_1(u) \ge f'_2(u)$.
If $u \not\in \min (\supp(f_1))$ and $u \not\in \min (\supp(f_2))$, 
then we have $f'_1 = f_1$ and $f'_2 = f_2$, thus $f'_1(u) \ge f'_2(u)$.

Next we prove (2).
If $u \in \supp(f'_2)$ and $u$ is minimal in $\supp(f'_1)$, 
then it follows from Lemma~\ref{lem:minimal} that $u$ is minimal in $\supp(f'_{2})$, 
hence we see that $f'_1(u) = - f_1(u) = - f_2(u) = f'_2(u)$.
This completes the proof.
\qed
\end{demo}

Now we are in position to prove Theorem~\ref{thm:e-triangulation} (b).

\begin{demo}{Proof of Theorem~\ref{thm:e-triangulation} (b)}
We need to show the following four claims:
\begin{enumerate}
\item[(1)]
If $K$ is a chain in $\eFF(P)$, then $T^{(e)}_K$ is a unimodular simplex.
\item[(2)]
If $K$ is a chain in $\eFF(P)$, then $T^{(e)}_K \subset \eCC(P)$.
\item[(3)]
$\bigcup_K T^{(e)}_K = \eCC(P)$, where $K$ runs over all chains in $\eFF(P)$.
\item[(4)]
If $K$ and $L$ are chains in $\eFF(P)$, then 
$T^{(e)}_K \cap T^{(e)}_L = T^{(e)}_{K \cap L}$.
\end{enumerate}

Recall that $T^{(e)}_K = \conv \ePhi(K)$ 
and $R_\varphi : \Real^P \to \Real^P$ is a linear map given by 
\[
\left( R_\varphi g \right)(v)
 = 
\begin{cases}
 g(v) &\text{if $\varphi(v) = 1$ or $0$,} \\
 -g(v) &\text{if $\varphi(v) = -1$,}
\end{cases}
\]
for $\varphi : P \to \{ 1, 0, -1 \}$.

(1)
If we put $C = \supp(K)$ and $\varphi = \sgn(K)$, 
then $T^{(e)}_K = R_\varphi(T_C)$ by Proposition~\ref{prop:simplex-reflection}.
Since $T_C$ is a unimodular simplex (Theorem~\ref{thm:triangulation} (b)) 
and $R_\varphi$ is a composition of reflections, 
we see that $T^{(e)}_K$ is a unimodular simplex.

(2)
We put $C = \supp(K)$ and $\varphi = \sgn(K)$, 
and take $\ep \in \{ 1, -1 \}^P$ such that $\varphi|_{\supp(\varphi)} = \ep|_{\supp(\varphi)}$.
Then, by using Lemma~\ref{lem:3-3} (c) and Lemma~\ref{lem:3-1}, we have 
\[
T^{(e)}_K
 = 
R_\ep (T_C)
\subset
R_\ep (\CC (P))
 =
\eCC_\ep(P)
\subset
\eCC(P).
\]

(3)
By using Lemma~\ref{lem:3-1} and Theorem~\ref{thm:triangulation} (b), we have
\[
\eCC(P)
 =
\bigcup_{\ep \in \{ 1, -1 \}^P} \eCC_\ep(P)
 =
\bigcup_{\ep \in \{ 1, -1 \}^P} R_\ep ( \CC(P) )
 =
\bigcup_{\ep \in \{ 1, -1 \}^P} \bigcup_C R_\ep ( T_C ),
\]
where $C$ runs over all chain of order filters of $P$.
Given a chain $C$ of order filters of $P$ and $\ep \in \{ 1, -1 \}^P$, 
we define $\varphi : P \to \{ 1, 0, -1 \}$ by putting
\[
\varphi(v)
 =
\begin{cases}
 \ep(v) &\text{if $v$ is minimal in some $F_i$,} \\
 0 &\text{otherwise.}
\end{cases}
\]
Then it follows from Lemma~\ref{lem:3-3} (c) that $R_\ep T_C = T^{(e)}_K$, 
where $K$ is the chain in $\eFF(P)$ corresponding to $(C,\varphi)$ 
under the bijection of Proposition~\ref{prop:bijection}.

(4)
We put $C = \supp (K)$, $\varphi = \sgn (K)$, $D = \supp (L)$ and $\psi = \sgn(L)$.
Then we have $T^{(e)}_K \subset V_\varphi$ and $T^{(e)}_L \subset V_\psi$ 
by Lemma~\ref{lem:3-3} (b).
If we define $\eta : P \to \{ 1, 0, -1 \}$ by putting
\[
\eta(v)
 = 
\begin{cases}
 1 &\text{if $\varphi(v) = \psi(v) = 1$,} \\
 -1 &\text{if $\varphi(v) = \psi(v) = -1$,} \\
 0 &\text{otherwise,}
\end{cases}
\]
then we have $V_\varphi \cap V_\psi = V_\eta$.
Hence we have
\[
T^{(e)}_K \cap T^{(e)}_L = T^{(e)}_K \cap T^{(e)}_L \cap V_\eta.
\]
Since $T^{(e)}_K = \conv ( \ePhi(K) )$ by definition, 
and $V_\eta$ is a ``boundary'' of $V_\varphi$, we see that
\[
T^{(e)}_K \cap V_\eta
 =
\conv ( \ePhi(K) ) \cap V_\eta
 =
\conv ( \ePhi(K) \cap V_\eta ).
\]
We take $\ep \in \{ 1, -1 \}^P$ satisfying $\eta|_{\supp(\eta)} = \ep|_{\supp(\eta)}$.
Then we have $\ePhi(K) \cap V_\eta$, $\ePhi(L) \cap V_\eta \subset \eAC_\ep(P)$. 
Since $R_\ep$ gives a bijection between $\eAC_\ep(P)$ and $\AC(P)$, 
it follows from Lemma~\ref{lem:3-4} that there exists a chain $C'$ of order filters of $P$ 
such that $R_\ep(\ePhi(\chi_{C'})) = \ePhi(K) \cap V_\eta$.
Hence we have
\[
T^{(e)}_K \cap V_\eta
 =
\conv ( R_\ep ( \ePhi(\chi_{C'}) ) )
 =
R_\ep \conv ( \ePhi(\chi_{C'}) ).
\]
Similarly there exists a chain $D'$ of order filters of $P$ such that 
\[
T^{(e)}_L \cap V_\eta
 =
\conv ( R_\ep ( \ePhi(\chi_{D'}) ) )
 =
R_\ep \conv ( \ePhi(\chi_{D'}) ).
\]
Therefore we have
\begin{align*}
T^{(e)}_K \cap T^{(e)}_L
 &=
(T^{(e)}_K \cap V_\eta) \cap (T^{(e)}_L \cap V_\eta)
\\
 &=
R_\ep \conv (\ePhi(\chi_{C'})) \cap R_\ep \conv (\ePhi(\chi_{D'}))
\\
 &=
R_\ep \left( \conv (\ePhi(\chi_{C'})) \cap \conv (\ePhi(\chi_{D'})) \right)
\\
 &=
R_\ep ( T_{C'} \cap T_{D'} ).
\end{align*}
By Theorem~\ref{thm:triangulation} (b), we see that 
$T_{C'} \cap T_{D'} = T_{C' \cap D'} = \conv ( \ePhi(\chi_{C' \cap D'}) )$.
Hence we have
\begin{align*}
T^{(e)}_K \cap T^{(e)}_L
 &=
R_\ep \left( \conv ( \ePhi(\chi_{C' \cap D'}) ) \right)
\\
 &=
R_\ep \left( \conv ( \ePhi(\chi_{C'}) \cap \ePhi(\chi_{D'}) ) \right)
\\
 &=
\conv \left( R_\ep (\ePhi(\chi_{C'})) \cap R_\ep (\ePhi(\chi_{D'}) \right)
\\
 &=
\conv \left( (\ePhi(K) \cap V_\eta) \cap (\ePhi(L) \cap V_\eta) \right)
\\
 &=
\conv \left( \ePhi(K) \cap \ePhi(L) \cap V_\eta \right)
\\
 &=
\conv \left( \ePhi(K) \cap \ePhi(L) \right)
\\
 &=
\conv \left( \ePhi(K \cap L) \right)
\\
 &=
T^{(e)}_{K \cap L}.
\end{align*}
This completes the proof of Theorem~\ref{thm:e-triangulation} (b).
\qed
\end{demo}

We conclude this subsection with giving a set of defining inequalities of a facet $T^{(e)}_K$, 
where $K$ is a maximal chain in $\eFF(P)$.
Recall the result of Stanley \cite{Stanley1986} on the defining inequalities of facets 
of the triangulations of $\OO(P)$ and $\CC(P)$.
To a maximal chain $C = \{ F_0 \supsetneq F_1 \supsetneq \dots \supsetneq F_d \}$ 
of order filters of $P$, 
we associate a linear extension $(v_1, \dots, v_d)$ and chains $C_1, \dots, C_d$ of $P$ 
as follows.
The linear extension $(v_1, \dots, v_d)$ is defined by
\[
F_i = F_{i-1} \cup \{ v_i \}
\quad(i=1, \dots, d).
\]
The chain $C_i$ is given inductively by
\begin{enumerate}
\item[(i)]
If $v_i$ is minimal, then we put $C_i = \{ v_i \}$;
\item[(ii)]
If $v_i$ is not minimal and $j$ is the largest index satisfying $v_j \lessdot v_i$, 
then we put $C_i = \{ v_i \} \cup C_j$.
\end{enumerate}

\begin{prop}
\label{prop:facet}
(Stanley \cite[Section~5]{Stanley1986})
Let $C$ be a maximal chain of order filters of $P$.
Let $(v_1, \dots, v_d)$ be the associated linear extension of $P$ 
and $C_1, \dots, C_d$ the associated chains of $P$.
Then we have
\begin{enumerate}
\item[(a)]
The facet $S_C$ of the triangulation $\ST_P$ of $\OO(P)$ is given by
\[
S_C = \{ f \in \Real^P : 0 \le f(v_1) \le f(v_2) \le \dots \le f(v_d) \le 1 \}.
\]
\item[(b)]
If $f \in S_C$, then we have
\[
( \Phi(f) )(v_i) = f(v_i) - f(v_j),
\]
where $j$ is the largest index satisfying $v_j \lessdot v_i$.
\item[(c)]
If we define
\[
L^C_i (g) = \sum_{v \in C_i} g(v),
\]
then the facet $T_C$ of the triangulation $\TT_P$ of $\CC(P)$ is given by
\[
T_C = \{ g \in \Real^P : 0 \le L^C_1(g) \le L^C_2(g) \le \dots \le L^C_d(g) \le 1 \}.
\]
\item[(d)]
If $g \in T_C$, then we have
\[
( \Psi(g) ) (v_i) = \sum_{v \in C_i} g(v),
\]
where $\Psi : \CC(P) \to \OO(P)$ is the inverse transfer map.
\end{enumerate}
\end{prop}

\begin{corollary}
\label{cor:e-facet-C}
Let $K$ be a maximal chain in $\eFF(P)$ 
and put $C = \supp(K)$, $\ep = \sgn(K)$.
Let $C_1, \dots, C_d$ be the chains of $P$ associated to $C$, 
and define
\[
\tilde{L}^K_i(g) = \sum_{v \in C_i} \ep(v) g(v) \quad(g \in \Real^P).
\]
Then the face $T^{(e)}_K$ of the triangulation $\eTT_P$ of $\eCC(P)$ is given by
\begin{equation}
\label{eq:e-facet-C}
T^{(e)}_K
 =
\{ g \in \Real^P : 0 \le L^K_1(g) \le L^K_2(g) \le \dots \le L^K_d(g) \le 1 \}.
\end{equation}
\end{corollary}

\begin{demo}{Proof}
It follows from Proposition~\ref{prop:simplex-reflection} and Proposition~\ref{prop:facet} (c).
\qed
\end{demo}

\subsection{%
Triangulation of $\eOO(P)$.
}

In this subsection, we transfer the triangulation of $\eCC(P)$ to $\eOO(P)$ 
via the inverse map $\ePsi$ of the enriched transfer map $\ePhi$.
In order to prove Theorem~\ref{thm:e-triangulation} (a), 
it is enough to show that $S^{(e)}_K =\conv K= \ePsi (T^{(e)}_K)$ 
and it is a unimodular simplex.

\begin{lemma}
\label{lem:3-5}
Let $K$ be a maximal chain in $\eFF(P)$ 
and put $C = \supp(K)$, $\ep = \sgn(K)$.
Let $(v_1, \dots, v_d)$ be the linear extension and $C_1, \dots, C_d$ the chains of $P$ associated to $C$.
For $g \in T^{(e)}_K$, we have
\[
( \ePsi(g) )(v_i)
 =
g(v_i) + \sum_{v \in C_i \setminus \{ v_i \}} \ep(v) g(v).
\]
\end{lemma}

\begin{demo}{Proof}
Since $T^{(e)}_K \subset V_\ep$ by Lemma~\ref{lem:3-3} (b), 
we have $|g(v)| = \ep(v) g(v)$ for $g \in T^{(e)}_K$ and $v \in P$, 
thus $|g| \in T_C$.
By Proposition~\ref{prop:facet}, we see that
\[
\max \{ S(|g|;B) : B \in \MC(P_{\le v_j}) \}
 =
\sum_{v \in C_j} |g(v)|
 =
\sum_{v \in C_j} \ep(v) g(v),
\]
where we recall $S(f;C) =
\sum_{v\in C}|f(v)|$ for $f \in \Real^P$ and a chain $C$ of $P$, 
and $\MC(P_{\le v_j})$ is the set of maximal chains of the subposet $P_{\le v_j}$.
Hence we obtain the desired identity.
\qed
\end{demo}

\begin{demo}{Proof of Theorem~\ref{thm:e-triangulation} (a)}
If $K$ is a maximal chain in $\eFF(P)$, then it follows from Lemma~\ref{lem:3-5} 
that $\ePsi$ is a unimodular linear map on $T^{(e)}_K$.
Hence, if $L$ is a chain in $\eFF(P)$ contained in $K$, 
then we see that $S^{(e)}_L = \ePsi(T^{(e)}_L)$ is a unimodular simplex 
because $T^{(e)}_L$ is a unimodular simplex (Theorem~\ref{thm:e-triangulation} (b)).
\qed
\end{demo}

We can use Lemma~\ref{lem:3-5} to give a set of defining inequalities of a facet $S^{(e)}_K$, 
where $K$ is a maximal chain in $\eFF(P)$.

\begin{prop}
\label{prop:e-facet-O}
Let $K$ be a maximal chain in $\eFF(P)$ 
and put $C = \supp(K)$, $\ep = \sgn(K)$.
Let $(v_1, \dots, v_d)$ be the linear extension and $C_1, \dots, C_d$ the chains of $P$ associated to $C$.
If we put
\[
\tilde{M}^K_i(f)
 =
\sum_{l=1}^r \ep(u_l) \prod_{j=l+1}^r (1 - \ep(u_j)) f(u_l)
\quad(f \in \Real^P),
\]
where $C_i = \{ u_1 \lessdot u_2 \lessdot \dots \lessdot u_r = v_i \}$, 
then the face $S^{(e)}_K$ of the triangulation $\eST_P$ of $\eOO(P)$ is given by
\begin{equation}
\label{eq:e-facet-O}
S^{(e)}_K
 =
\{ f \in \Real^P : 0 \le \tilde{M}^K_1(f) \le \tilde{M}^K_2(f) \le \dots \le \tilde{M}^K_d(f) \le 1 \}.
\end{equation}
\end{prop}

\begin{demo}{Proof}
It is easy to prove by induction on $k$ that
\[
f(u_k) = g(u_k) + \sum_{i=1}^k \ep(u_i) g(u_i) \quad(k=1, \dots, r)
\]
if and only if
\[
g(u_k) = f(u_k) - \sum_{i=1}^{k-1} \ep(u_i) \prod_{j=i+1}^{k-1} (1 - \ep(u_j)) f(u_i)
\quad(k=1, \dots, r).
\]
Hence we have
\[
\tilde{M}^K_i(f) = \tilde{L}^K_i(\ePhi(f)).
\]
On the other hand, by Theorem~\ref{thm:e-triangulation} (b) and Corollary~\ref{cor:e-facet-C}, 
we see that $f \in S^{(e)}_K$ if and only if 
\[
0 \le \tilde{L}^K_1(\ePhi(f)) \le \cdots \le \tilde{L}^K_d(\ePhi(f)) \le 1.
\]
Hence we obtain (\ref{eq:e-facet-O}).
\qed
\end{demo}

\section{%
Identification with Ohsugi--Tsuchiya's triangulations
}
\label{sec:algebraic}

Ohsugi--Tsuchiya \cite{OT1,OT2} computed the squarefree initial ideals of  
the toric ideals of $\eOO(P)$ and $\eCC(P)$ 
with respect to certain monomial orderings. 
This gives regular unimodular triangulations of $\eOO(P)$ and $\eCC(P)$.
In this section, we show that these triangulations coincide with the triangulations 
given in Theorem~\ref{thm:e-triangulation}.


Let $\K[\vectx] = \K[x_1,\dots, x_n]$ be the polynomial ring in $n$ variables 
$x_1, \dots, x_n$ over a field $\K$
and $\Delta$ a simplicial complex on $[n]:=\{1,2\ldots,n\}$.
To a subset $F \subset [n]$, we associate a monomial
\[
\vectx_F=\prod_{i \in F} x_i.
\]
The \emph{Stanley--Reisner ideal} of $\Delta$ is the ideal $I_\Delta$ of $\K[\vectx]$ which is generated by those squarefree monomials $\vectx_F$ with $F \notin \Delta$.
On the other hand, given an arbitrary squarefree monomial ideal $I$ of $\K[\vectx]$, there is a unique simplicial complex $\Delta(I)$ such that $I=I_{\Delta(I)}$.

\subsection{%
Triangulation of enriched order polytope
}

In this subsection, we prove that the triangulation $\eST_P$ of $\eOO(P)$ given 
in Theorem~\ref{thm:e-triangulation} (a) 
coincides with the one algebraically defined in \cite{OT2}.

Let $P$ be a finite poset with $d$ elements and let $R[\eOO] = \K [ \{ x_f : f \in \eFF(P) \} ]$ be the polynomial ring 
in the variables $x_f$ ($f \in \eFF(P)$).

\begin{prop}
\label{prop:init_ideal_eO}
(\cite[Theorem 5.2]{OT2})
Let $I_{\eOO(P)}$ be the ideal of $R[\eOO]$ generated by all squarefree monomials
$x_f x_g$ satisfying either of the following conditions:
\begin{enumerate}
\item[(i)]
there exists $v \in \min(\supp(f)) \cap \min(\supp(g))$ such that $f(v) \neq g(v)$;
\item[(ii)]
$\supp(f) \not\sim \supp(g)$ and $f(v) = g(v)$ 
for each $v \in \min(\supp(f))\cap \min(\supp(g))$,
\end{enumerate}
where the symbol $A \not\sim B$ means that $A \nsubseteq B$ and $A \nsupseteq B$. 
Then $\Delta(I_{\eOO(P)})$ is a regular unimodular triangulation of $\eOO(P)$.
\end{prop}

Now we can show that this triangulation $\Delta(I_{\eOO(P)})$ 
coincides with the triangulation given in Theorem~\ref{thm:e-triangulation} (a).

\begin{prop}
\label{prop:sametriangulation_eO}
With the notations above, we have
$\Delta(I_{\eOO(P)})=\eST_P$.
\end{prop}

\begin{demo}{Proof}
Since both of $\Delta(I_{\eOO(P)})$ and 
$\eST_P$
are unimodular triangulations of $\eOO(P)$, 
the number of maximal simplices are same.
Hence it is enough to show that $x_{f_0} \cdots x_{f_d} \not\in I_{\eOO(P)} $ 
for any maximal chain $K = \{ f_0 \gtrdot \dots \gtrdot f_d \}$ of $\eFF(P)$.

Let $K=\{ f_0 > \cdots > f_d\}$ be a maximal chain of $\eFF(P)$, 
and assume to the contrary that $x_{f_0} \cdots x_{f_d} \in I_{\eOO(P)}$.
Then there exists a pair of indices $i<j$ such that $f_i$ and $f_j$ satisfy 
the condition (i) or (ii) in Proposition~\ref{prop:init_ideal_eO}.
Since $f_i > f_j$ in $\eFF(P)$, one has $\supp(f_i) \supsetneq \supp(f_j)$ 
and $f_i$ and $f_j$ do not satisfy the condition (ii).
Hence there exists $v \in \min(\supp(f_i)) \cap \min(\supp(f_j))$ with $f_i(v) \neq f_j(v)$.
However, since $f_i > f_j$ in $\eFF(P)$ and $v \in \supp(f_j)$ and $v$ is minimal in $\supp(f_i)$, 
we obtain $f_i(v) = f_j(v)$, which is a contradiction.
Thus it follows that $x_{f_0} \cdots x_{f_d} \not\in I_{\eOO(P)}$.
\qed
\end{demo}



\subsection{%
Triangulation of enriched chain polytope
}

In this subsection, we prove that the triangulation $\eTT_P$ of $\eCC(P)$ 
given in Theorem~\ref{thm:e-triangulation} (b) 
coincides with the one algebraically defined in \cite{OT1}.
Let $R[\eCC]$ be the polynomial ring in variables $y_g$ ($g \in \eAC(P)$).

\begin{prop}
\label{prop:init_ideal_eC}
(\cite[Theorem 1.4]{OT1})
Let $I_{\eCC(P)}$ be the ideal of $R[\eCC]$ generated 
by all squarefree monomials $y_g y_h$ satisfying either of the following conditions:
\begin{enumerate}
\item[(i)]
there exists $v \in \supp(g) \cap \supp(h)$ such that $g(v) \neq h(v)$;
\item[(ii)]
$\langle \supp(g) \rangle \not\sim \langle \supp(h) \rangle$ 
and $g(v) = h(v)$ for any $v \in \supp(g) \cap \supp(h)$.
\end{enumerate}
Then  $\Delta ( I_{\eCC(P)} )$ is a regular unimodular triangulation of $\eCC(P)$.
\end{prop}

Finally, we show that this triangulation $\Delta(I_{\eCC(P)})$ 
coincides with the triangulation given in Theorem~\ref{thm:e-triangulation} (b).

\begin{prop}
\label{prop:sametriangulation_eC}
With the same notation as above, we have
$\Delta (I_{\eCC(P)}) =\eTT_P$.
\end{prop}

\begin{demo}{Proof}
This follows from the fact that the map $x_f \mapsto y_{\ePhi(f)}$ induces the ring isomorphism 
\[
\dfrac{R[\eOO(P)]}{I_{\eOO(P)}} \cong \dfrac{R[\eCC(P)]}{I_{\eCC(P)}}.
\]
\qed
\end{demo}



\subsubsection*{Declarations}

\subsubsection*{Conflict of Interest}
On behalf of all authors, the corresponding author states that there is no conflict of interest.

\end{document}